\documentclass[11pt]{article}
\usepackage{amsmath}
\usepackage{amsfonts}
\usepackage{amsbsy}
\usepackage{amsthm}
\usepackage{amssymb}
\usepackage{graphicx}
\usepackage{psfrag}
\usepackage{epsfig}
\usepackage{citesort}
\usepackage{graphicx}

\newtheorem{theorem}{Theorem}

\newtheorem{definition}{Definition}

\theoremstyle{remark}

\begin{document}

\title{Second-kind integral solvers for TE and TM problems of
  diffraction by open arcs}

\author{Oscar P. Bruno and St\'{e}phane K. Lintner}
        
\maketitle
\begin{abstract}
  We present a novel approach for the numerical solution of problems of
  diffraction by open arcs in two dimensional space.  Our methodology relies
  on composition of {\em weighted versions} of the classical integral
  operators associated with the Dirichlet and Neumann problems (TE and TM
  polarizations, respectively) together with a generalization to the
  open-arc case of the well known closed-surface Calder\'on formulae.  When
  used in conjunction with spectrally accurate discretization rules and
  Krylov-subspace linear algebra solvers such as GMRES, the new second-kind
  TE and TM formulations for open arcs produce results of high accuracy in
  small numbers of iterations---for low and high frequencies alike.
\end{abstract}

\section{Introduction}
Problems of diffraction by {\em infinitely thin open surfaces} play
central roles in a wide range of problems in science and engineering,
with important applications to antenna and radar design, electronics,
optics, etc. Like other wave scattering problems, open surface
problems can be treated by means of numerical methods that rely on
approximation of Maxwell's Equations over volumetric domains (on the
basis of, e.g., finite-difference or finite-element methods) as well
as methods based on boundary integral equations. As a result of the
singular character of the electromagnetic fields in the vicinity of
open edges, open-surface scattering configurations present major
difficulties for both volumetric and boundary integral methods.
Boundary integral approaches, which require discretization of domains
of lower dimensionality than those involved in volumetric methods, can
generally be treated efficiently, even for high-frequencies, by means
of accelerated iterative scattering
solvers~\cite{BleszynskiBleszynskiJaroszewicz,BrunoKunyansky,Rokhlin}.
Unfortunately, integral methods for open surfaces do not give rise,  at least
in their classical formulations, to Fredholm integral operators of the second-kind, and can therefore prove to be computationally expensive---as
the eigenvalues of the resulting equations accumulate at zero and/or
infinity and, thus, iterative solution of these equations often
requires large numbers of iterations.

This paper presents new Fredholm integral equations of {\em second-kind} and
associated numerical algorithms for problems of scattering by
two-dimensional open arcs $\Gamma$ (i.e., infinite cylinders of
cross-section $\Gamma$) under either transverse electric (TE) or transverse
magnetic (TM) incident fields.  The new second-kind Fredholm equations
(which result from composition of appropriately modified versions
$\mathbf{S}_\omega$ and $\mathbf{N}_\omega$ of the classical single-layer
and hypersingular integral operators $\mathbf{S}$ and $\mathbf{N}$) provide,
for the first time, a generalization of the classical closed-surface
Calder\'on formulas to the open-arc case. In particular, the new
formulations possess highly favorable spectral properties: their eigenvalues
are highly clustered, and they remain bounded away from zero and infinity,
even for problems of very high frequency. When used in conjunction with
spectrally accurate discretization rules and Krylov-subspace linear algebra
solvers such as GMRES, the new open-arc formulations produce results of high
accuracy in small numbers of iterations---for low and high frequencies
alike.

The new second-kind formulation for the TM problem is particularly
beneficial, as it gives rise to order-of-magnitude improvements in computing
times over the corresponding weighted hypersingular formulation. Such gains
do not occur in the TE case: although the new second-kind TE equation
requires fewer iterations than the corresponding weighted first kind
formulation, the total computational cost of the second-kind equation is
generally higher in the TE case---since the application of the first-kind
operator can be significantly less expensive than the application of the
composite second-kind operator.

The difficulties that arise as integral equations are used to treat open
surface scattering problems are of course well known, and many contributions
have been devoted to their treatment. Like the present work,
references~\cite{PovznerSuharesvki} and~\cite{ChristiansenNedelec} seek to
tackle these problems by means of generalizations of the classical
Calder\'on relations to the case of open surfaces. The first of these
contributions establishes that the combination $\mathbf{N}\mathbf{S}$ can be
expressed in the form $\mathbf{I} + \mathbf{T}_K$, where the kernel
$\mathbf{K}(x,y)$ of the operator $\mathbf{T}_K$ has local singularity of at
most $O\left(\frac{1}{|x-y|}\right)$. This early result however does not
take into account the singular edge behavior; the resulting operator
$\mathbf{T}_K$ is not compact (in fact it gives rise to extreme
singularities at the edge~\cite{LintnerBruno}) and $\mathbf{I} +
\mathbf{T}_K$ is therefore not a second-kind operator in any numerically
meaningful functional space. When used in conjunction with boundary elements
that vanish on the edges (by means of well-chosen projections) however, the
combination $\mathbf{N}\mathbf{S}$ does give rise to reduction of iteration
numbers, as demonstrated in reference~\cite{ChristiansenNedelec} through
numerical examples for low frequency problems. This contribution does not
include details on accuracy, and it does not utilize integral weights to
resolve the solution's edge singularity.  A related but different method was
introduced in~\cite{AntoineBendaliDarbasMarion} which exhibits, once again,
low iteration numbers for low-frequency problems, but which does not resolve
the singular edge behavior and for which no accuracy studies have been
presented. Finally, high-order integration rules for the single-layer and
hypersingular operators adapted to open arcs were introduced in a Galerkin
framework in~\cite{StephanWendland,StephanHsiao,StephanTran}. These methods
have thus far only been applied for simple geometries and at low
frequencies, and limited information is available on the actual convergence
properties and performance of their computational implementations.

A second class of methods include those proposed in
references~\cite{AtkinsonSloan,Monch,RokhlinJiang}. The
contribution~\cite{AtkinsonSloan}, some aspects of which are incorporated in
our method, treats the Dirichlet problem for Laplace's equation by means of
second kind equations. The basis of this approach lies in the observation
that the cosine basis has the dual positive effect of diagonalizing the
logarithmic potential for a straight arc and removing the singular edge
behavior---so that the inverse of the logarithmic potential can be easily
computed and used as a preconditioner to produce a second kind operator for
a general arc.  The approach~\cite{Monch}, which also uses a cosine basis,
treats the Neumann problem for the non-zero frequency Helmholtz equation
with spectral accuracy by means of first kind equations.  The
contribution~\cite{RokhlinJiang}, finally, treats, just
like~\cite{AtkinsonSloan}, the Laplace problem by means of second-kind
equations resulting from inversion of the straight arc logarithmic
potential; like~\cite{Monch}, further, it produces spectral accuracy through
use of the cosine transforms.  The second-kind integral approach developed
in~\cite{AtkinsonSloan} and later revisited in \cite{RokhlinJiang} seems
essentially limited to the specific problem for which it was proposed:
neither an extension to the Neumann problem nor to the full three
dimensional problem seem straightforward. And, more importantly, this
approach does not lead to adequately preconditioned equations for non-zero
frequencies: a simple experiment conducted in Section~\ref{sec_eigenvalues}
shows that a direct generalization of the algorithm~\cite{AtkinsonSloan} to
the Helmholtz problem generally requires significantly {\em more} linear
algebra iterations than are necessary if the operator $\mathbf{S}_\omega$
alone is used.


The remainder of this paper is organized as follows: after recalling in
Section~\ref{sec_preliminaries} the classical boundary integral formulations
for the TE and TM open-arc scattering problems, in Sections~\ref{A}
and~\ref{B} we present our new weighted operators $\mathbf{S}_\omega$ and
$\mathbf{N}_\omega$ as well as certain periodized counterparts $\tilde{S}$
and $\tilde{N}$ (which are obtained by considering a sinusoidal changes of
variables for source and observation points). Our main result, the
generalized Calder\'on formula, is presented in
Section~\ref{C}. Section~\ref{D} then presents results of numerical
evaluation of eigenvalues for a non-trivial open arc problem, illustrating
the spectral properties of previous open-arc operators as well as the second
kind operators introduced in this paper.  Theoretical considerations
concerning the new second-kind equations are presented in
Section~\ref{sec_theory}, including a succinct but complete proof of the
open-arc Calder\'on formulae; a more detailed theoretical discussion,
including full mathematical technicalities, can be found
in~\cite{LintnerBruno}. The high-order quadrature rules we use for
evaluation of the new integral operators are described in
Section~\ref{sec_quadrature}. Numerical results, finally, are presented in
Section~\ref{sec_results} for a wide range of frequencies and for various
geometries (including a brief study of resonant open cavities) demonstrating
the uniformly well conditioned character of the integral formulations
proposed in this paper.

\section{Preliminaries\label{sec_preliminaries}}
Let $\Gamma$ denote a smooth open arc in the plane.  The TE and TM problems
of scattering by the open arc $\Gamma$ amount to two-dimensional problems
for the Helmholtz equation with Dirichlet and Neumann boundary conditions on
$\Gamma$, respectively:
\begin{equation}\label{b_conds}
  \left\{ \begin{array}{llll} \Delta u +k^2 u = 0 \quad\mbox{outside}\quad
      \Gamma , & u|_{\Gamma} = f & \mbox{(TE)}\\
      \Delta v +k^2 v = 0 \quad\mbox{outside}\quad \Gamma, & \frac{\partial
        v}{\partial \textbf{n}}|_{\Gamma} = g & \mbox{(TM)}.
\end{array}\right.
\end{equation}
Here $\textbf{n}$ denotes the normal to $\Gamma$, and $f$ and $g$ are given
in terms of the incident electric excitation $u^{inc}$: $f=-u^{inc}$ and
$g=-\partial u^{inc}/\partial \textbf{n}$ on $\Gamma$.

\subsection{Classical open-arc equations for TE and TM problems}
The (unique) solutions $u$ and $v$ of the TE and TM problems above can be
expressed, for $\mathbf{r}\not\in\Gamma$, in terms of single- and
double-layer potentials of the form $u(\mathbf{r})= \int_\Gamma
G_k(\mathbf{r},\mathbf{r}')\mu(\mathbf{r}') d\ell'$ and $v(\mathbf{r})=
\int_\Gamma \frac{\partial G_k(\mathbf{r},\mathbf{r}')}{\partial
\textbf{n}_{\mathbf{r}'}}\nu(\mathbf{r}') d\ell'$; see
e.g.~\cite{StephanWendland,StephanWendland2,VainikkoSaranen}.
 Here $\textbf{n}_{\mathbf{r}'}$ is a unit vector normal to $\Gamma$
 at the point $\mathbf{r}'\in\Gamma$ (we assume, as we may, that
 $\textbf{n}_{\mathbf{r}'}$ is a smooth function of
 $\mathbf{r}'\in\Gamma$), and, calling $H_0^1$ the first Hankel
 function of order zero,
\begin{equation}\label{Gk_def}
G_k(\mathbf{r},\mathbf{r}')=\left\{\begin{array}{ll} \frac{i}{4}H_0^1(k
|\mathbf{r}-\mathbf{r}'|),& k > 0 \\
-\frac{1}{2\pi}\ln|\mathbf{r}-\mathbf{r}'|,&k=0 \end{array}\right. 
\end{equation}
and $\frac{\partial G_k(\mathbf{r},\mathbf{r}')}{\partial
  \mathbf{n}_{\mathbf{r}'}}=\mathbf{n}_{\mathbf{r}'}\cdot\nabla_{\mathbf{r}'}G_k(\mathbf{r},\mathbf{r}').$
Denoting by $\mathbf{S}$ and $\mathbf{N}$ the single-layer and
hypersingular operators
\begin{equation}\label{Sdef}
 \mathbf{S}[\mu](\mathbf{r})=  \int_\Gamma
 G_k(\mathbf{r},\mathbf{r}')\mu(\mathbf{r}') d\ell'\quad , \quad
 \mathbf{r}\in\Gamma,
 \end{equation}
and
\begin{equation}\label{Ndef}
\begin{split}
  \mathbf{N}[\nu](\mathbf{r})= &\; \frac{\partial }{\partial
  \textbf{n}_{\mathbf{r}}}\int_\Gamma \frac{\partial
  G_k(\mathbf{r},\mathbf{r}')}{\partial
  \textbf{n}_{\mathbf{r}'}}\nu(\mathbf{r}') d\ell'\\
  \stackrel{\mathrm{def}}{=}&\lim\limits_{z\rightarrow 0^+} \frac{\partial
  }{\partial z}\int_\Gamma \frac{\partial
  G_k(\mathbf{r}+z\mathbf{n}_{\mathbf{r}},\mathbf{r}')}{\partial
  \textbf{n}_{\mathbf{r}'}}\nu(\mathbf{r}') d\ell'\quad , \quad
  \mathbf{r}\in\Gamma,
\end{split}
\end{equation}
further, the densities $\mu$ and $\nu$ are the unique solutions of the
integral equations
\begin{equation}\label{bad}
\mathbf{S}[\mu]=f, \quad
\mathbf{N}[\nu]=g.
\end{equation}

It is known that the eigenvalues of the integral operators in
Equations~\eqref{bad} accumulate at zero and infinity, respectively (see
e.g. Figure~\ref{EigsFigSpirale}). As a result (and as illustrated in
Section~\ref{sec_results}), solution of these equations by means of
Krylov-subspace iterative solvers such as GMRES generally require large
numbers of iterations. In addition, as discussed in Section~\ref{edge_beh},
the solutions $\mu$ and $\nu$ of equations~\eqref{bad} are not smooth at the
end-points of $\Gamma$, and, thus, they give rise to low order convergence
(and require high discretization of the densities for a given accuracy) if
standard quadrature methods are used.  Details in these regards are provided
in what follows.

\subsection{Calder\'on relation for closed and open arcs.} 
As is well known~\cite{ColtonKress1,Nedelec,ColtonKress2} second-kind
formulations for closed surfaces (arcs) result from: (1)~The classical
jump relations associated with singular integral operators, as well as
(2)~application of the Calder\'on formula
\begin{equation}\label{CalderonClosed}
\mathbf{N}_c\mathbf{S}_c=-\frac{\mathbf{I}}{4}+\mathbf{K}_c,
\end{equation}
which shows that the composition of the closed-surface hypersingular
and single-layer operators $\mathbf{N}_c$ and $\mathbf{S}_c$ equals a
compact perturbation $\mathbf{K}_c$ of the identity.

In the open-surface context however, jump relations cannot be
exploited since, according to the boundary conditions~\eqref{b_conds},
the same limits must be reached on both sides of the open surface
$\Gamma$. Further, use of the combined operator $\mathbf{N}\mathbf{S}$
does not lead to well-posed equations: as shown
in~\cite{LintnerBruno}, the image under $\mathbf{N}\mathbf{S}$ of the
constant function $1$ is highly singular, with edge asymptotics $
\mathbf{NS}[ 1 ](\mathbf{r})=O\left(\frac{1}{d(\mathbf{r})}\right)$,
$\mathbf{r}\in \Gamma$, where $d(\mathbf{r})$ denotes the Euclidian
distance to the edge. The combination $\mathbf{N}\mathbf{S}$ is also
problematic in the functional framework set forth
in~\cite{StephanWendland,StephanWendland2,Stephan}, since the image
under $\mathbf{S}$ of the Sobolev space $H^\frac{1}{2}(\Gamma)$ (as
defined in~\cite{Stephan}) is larger than the domain of definition of
$\tilde{H}^{\frac{1}{2}}(\Gamma)$ (as defined in~\cite{Stephan}) of
the operator $\mathbf{N}$.
\subsection{Regularity and singular behavior at the edge\label{edge_beh}}
The singular character of the solutions of equations~\eqref{bad} is well
documented~\cite{Stephan,Costabel}: $\mu$ and
$\nu$ can be expressed in the forms
\begin{equation}\label{StephanExp1}
\mu\sim \frac{\chi_1}{\sqrt{d}} +\eta, \quad
\nu\sim \chi_2\sqrt{d} +\zeta,
\end{equation}
where $d$ denotes the distance to the edge, $\chi_1$ and $\chi_2$ are
smooth cut-off functions, and where the functions $\eta$ and $\zeta$
are somewhat smoother than $\mu$ and $\nu$.  More recently
furthermore, it was established~\cite{Costabel} that expressions such
as~(\ref{StephanExp1}) can be viewed as part of an expansion in powers
of $d^{1/2}$ which can be carried out as long as the degree of
smoothness of the surface and the right hand side of the equations
allow it. For example, if the curve itself and the functions $f$ and
$g$ in~\eqref{bad} are infinitely differentiable, it follows that
\begin{equation}\label{CostabelExp1}
\mu=\frac{\alpha}{\sqrt{d}}, \quad
\nu=\beta\sqrt{d},
\end{equation}
where $\alpha$ and $\beta$ are infinitely differentiable functions throughout
$\Gamma$, \textit{up to and including the endpoints}~\cite{Costabel}. Thus
the singular character of these solutions is fully characterized by the
factors $d^{1/2}$ and $d^{-1/2}$ in equation~\eqref{CostabelExp1}.
\section{Well-posed Second-Kind Integral Equation Formulations}\label{sec_calderon}

\subsection{Weighted operators\label{A}}
In view of the regularity results~\eqref{CostabelExp1}, we introduce a
positive integral weight $\omega(\mathbf{r})>0$ with asymptotic edge
behavior $\omega\sim d^{1/2}$
at the edge (by which it is implied that the quotient $\omega /
d^{1/2}$ is infinitely differentiable up to the edge), and we define
the weighted operators
\begin{equation}\label{Somega}
\mathbf{S}_\omega[\alpha]=  \mathbf{S}\left[\frac{\alpha}{\omega}\right], \quad
\mathbf{N}_\omega[\beta] = \mathbf{N}[\omega\beta],
\end{equation}
so that, in view of the discussion of Section~\ref{edge_beh}, for smooth
functions $f$ and $g$, the solutions of the equations
\begin{equation}\label{good}
\mathbf{S}_\omega[\alpha]= f \quad \mbox{and}\quad
\mathbf{N}_\omega[\beta]= g
\end{equation}
are also smooth. 

Without loss of generality, we choose a smooth parameterization
$\mathbf{r}(t)=\left(x(t),y(t)\right)$ of $\Gamma$ defined in the interval
$[-1,1]$, and we define the weight $\omega$ canonically as
$\omega(\mathbf{r}(t))= \sqrt{1-t^2}.$
As a result, the operators $\mathbf{S}_\omega$ and $\mathbf{N}_\omega$ give rise to the 
parameter-space operators
\begin{equation}\label{sOmegaDef}
S_\omega[\varphi](t)=\int_{-1}^1 G_k\left(\mathbf{r}(t),\mathbf{r}(t')\right)\frac{
\varphi(t')}{\sqrt{1-t'^2}}\;\tau(t') dt',
\end{equation}
and
\begin{equation}\label{nOmegaDef}
\begin{split}
N_\omega[\psi](t)=\lim\limits_{z\rightarrow
0^+}\frac{\partial}{\partial z}\int_{-1}^1 
\frac{\partial}{\partial
\textbf{n}_{\mathbf{r}(t')}}G_k\left(\mathbf{r}(t)+z\textbf{n}_{\mathbf{r}(t)},\mathbf{r}(t')\right)\\
\psi (t') \tau(t')\sqrt{1-t'^2}dt',
\end{split}
\end{equation}
defined on functions $\varphi$ and $\psi$ of the variable $t$, $-1\leq t\leq
1$; in these equations we have set $\tau(t)=|\frac{d
\textbf{r}(t)}{dt}|$. Clearly, for $\varphi(t) = \alpha(\mathbf{r}(t))$ and
$\psi(t) = \beta(\mathbf{r}(t))$ we have
$\mathbf{S}_\omega[\alpha](\mathbf{r}(t)) = S_\omega[\varphi](t)$ and
$\mathbf{N}_\omega[\beta](\mathbf{r}(t)) = N_\omega[\psi](t).$

\subsection{Periodized spaces and regularity}\label{B}
Introducing the changes of variables $t=\cos\theta$ and
$t'=\cos\theta'$, and defining $\textbf{n}_\theta=
\textbf{n}_{\mathbf{r}(\cos\theta)}$, we obtain the periodic weighted
single-layer and hypersingular operators
\begin{equation}\label{ssdef}
  \tilde{S}[\gamma](\theta)=\int_{0}^\pi
  G_k(\mathbf{r}(\cos\theta),\mathbf{r}(\cos\theta'))\gamma(\theta')\tau(
  \cos\theta')d\theta'
\end{equation}
and 
\[
\tilde{N}[\gamma](\theta) =\lim\limits_{z\rightarrow
  0^+}\frac{\partial}{\partial z}\int_{0}^\pi \frac{\partial}{\partial
  \textbf{n}_{\theta'}}G_k(\mathbf{r}(\cos\theta)+z\textbf{n}_{\theta},\mathbf{r}(\cos\theta'))\\
\gamma(\theta') \tau(\cos\theta')\sin^2\theta'd\theta'.
\]
Clearly then, the solutions to the periodic equations
\begin{equation}\label{tildeEqs}
\tilde{S}[\tilde{\varphi}]= \tilde f, \quad 
\tilde{N}[\tilde{\psi}] = \tilde g, 
\end{equation}
where $\tilde f(\theta) = f(\mathbf{r}(\cos\theta))$ and $\tilde
g(\theta) = g(\mathbf{r}(\cos\theta))$,
are related to the solutions of equations~\eqref{good} by the relations
\begin{equation}\label{phi_psi_theta}
  \tilde \varphi(\theta) = \varphi(\cos\theta)\quad,\quad \tilde \psi(\theta)
  = \psi(\cos\theta).
\end{equation}
In view of~\eqref{phi_psi_theta}, the solutions to
equations~\eqref{tildeEqs} are smooth and periodic, and it is therefore
natural to study the properties of $\tilde{S}$ and $\tilde{N}$ in the the
Sobolev spaces $H^s_e(2\pi)$ of $2 \pi$ periodic and even functions cf.
\cite{YanSloan,BrunoHaslam}.
\begin{definition}
  Let $s\in \mathbb{R}$. The Sobolev space $H^s_e(2\pi)$ is defined as
  the set of $2\pi$-periodic functions of the form
  $v(\theta)=\frac{1}{2}a_0 + \sum\limits_{m=1}^\infty a_m \cos(
  m\theta )$
  for which the $s$-norm $\|v\|_s^2 = |a_0|^2 + 2\sum\limits_{m=1}^\infty m^{2s}|a_m|^2$
is finite. 
\end{definition}
Clearly the set $\{\cos(n\theta):n\in \mathbb{N} \}$ is a basis of the
Hilbert space $H^s_e(2\pi)$ for all $s$.

\begin{figure}
\begin{center}
\includegraphics[width=0.24\textwidth]{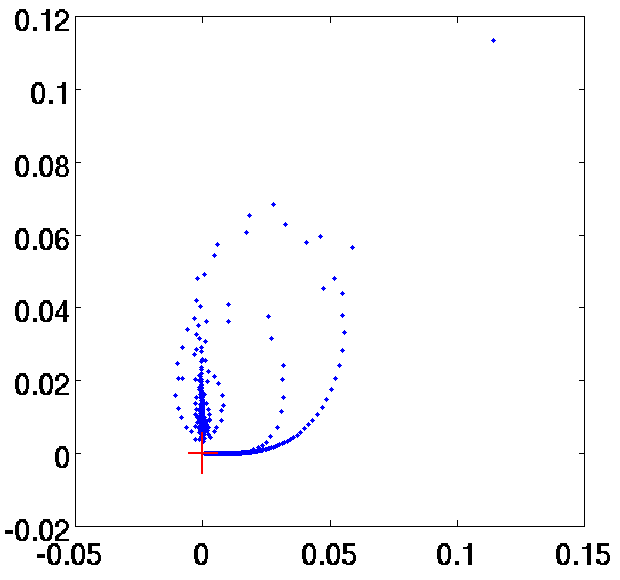}
\includegraphics[width=0.24\textwidth]{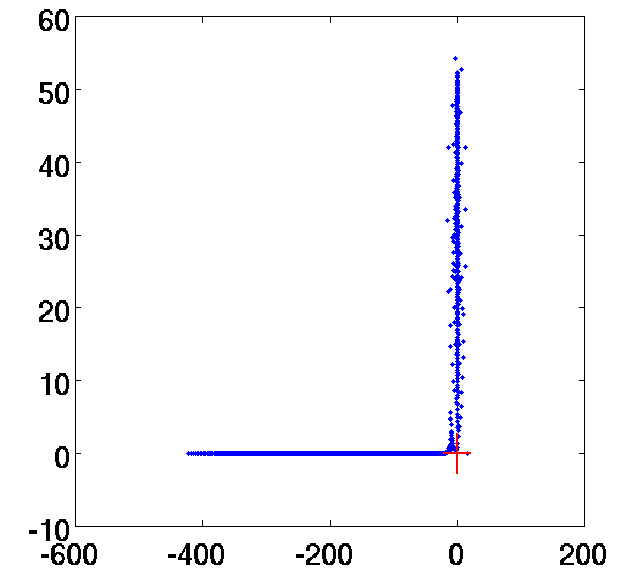}
\includegraphics[width=0.24\textwidth]{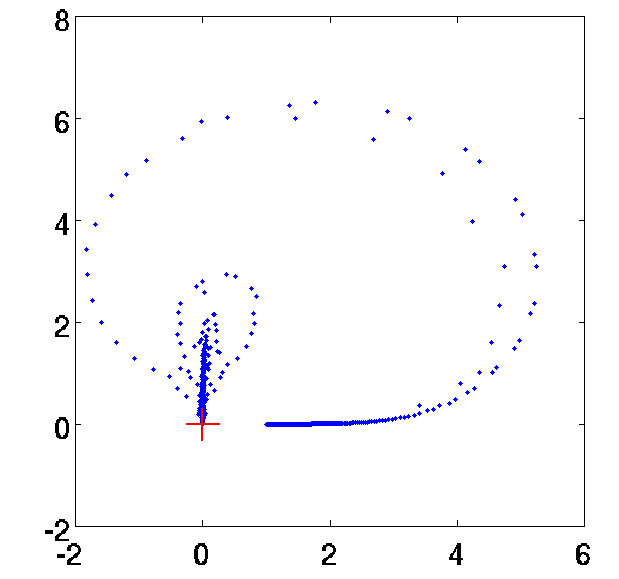}
\includegraphics[width=0.24\textwidth]{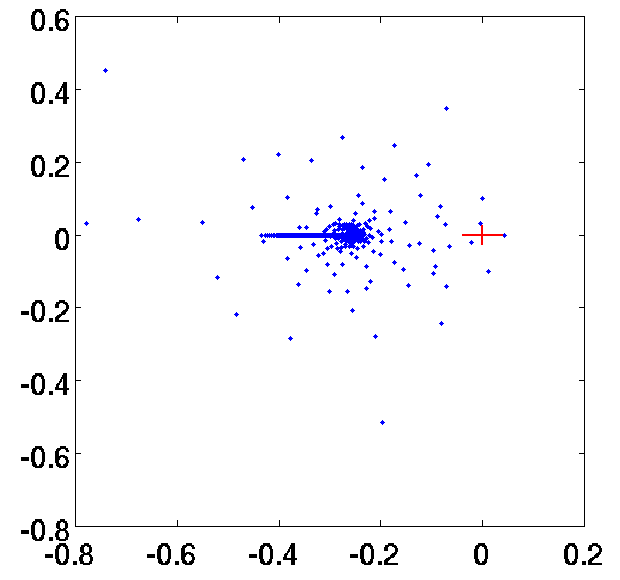}
\caption{Eigenvalue distribution for the spiral-shaped arc displayed in
  Figure~\ref{spiraleFigure}, with frequency $\frac{L}{\lambda}=100$, for
  the various operators under consideration.  Left $\tilde{S}$. Left-center
  $\tilde{N}$. Right-center $\tilde{S}_0^{-1}\tilde{S}$.  Right
  $\tilde{N}\tilde{S}$. Note the important difference of scale between the
  four plots.\label{EigsFigSpirale}}
\end{center}
\end{figure}

\subsection{Generalized Calder\'on Formula\label{C}}

Given the smoothness of the solutions of the equations arising from the
weighted periodic operators $~\tilde{S}$ and $~\tilde{N}$, and in light of
equation~\eqref{CalderonClosed}, it is reasonable to consider the
composition $\tilde{N}\tilde{S}$ as possible basis for solution of open-arc
problems. As stated in the following theorem and demonstrated below, this
combination does indeed give rise to a second-kind integral formulation. The
statement of the theorem makes use of the concept of bicontinuity: an
operator $L$ between two Hilbert spaces is said to be bicontinuous if and
only if $L$ is invertible and both $L$ and $L^{-1}$ are continuous
operators.

\begin{theorem}\label{NStheorem}\label{th_1}
  \
\begin{itemize} 
\item[(i)] For all $s>0$ the operators $\tilde{N}$ and $\tilde{S}$
  define bicontinuous mappings $\tilde{S}: H^s_e(2\pi) \rightarrow
  H^{s+1}_e(2\pi)$ and $\tilde{N}: H^{s+1}_e(2\pi) \rightarrow
  H^{s}_e(2\pi)$.
Thus, the composite operator $\tilde{N}\tilde{S}$ is bicontinuous
from $H^s_e(2\pi)$ into $H^s_e(2\pi)$ for all $s>0$.
\item[(ii)] The generalized Calder\'on formula
\begin{equation}\label{NSfact}
\tilde{N}\tilde{S}= \tilde J^\tau_0 +\tilde K
\end{equation}
holds, where ${\tilde K}: H^s_e(2\pi) \rightarrow H^s_e(2\pi)$ is a
compact operator, and where $\tilde J^\tau_0: H^s_e(2\pi) \rightarrow
H^s_e(2\pi)$ is a bicontinuous operator independent of the frequency
$k$.
\item[(iii)] The set $\sigma_s$ of eigenvalues of $\tilde{J}^\tau_0 :
  H^s_e(2\pi) \to H^s_e(2\pi)$ equals the union of the discrete set
  $\Lambda_\infty=\{\lambda_0=-\frac{\ln2}{4},\quad\lambda_n=-\frac{1}{4}-\frac{1}{4n}:n>0\}$
  and a certain open set $\Lambda_s$ which is bounded away from zero
  and infinity. The sets $\Lambda_s$ are nested, they form a
  decreasing sequence, ($\Lambda_s \varsubsetneq \Lambda_{s'}, s >
  s'$) and they satisfy
  $\bigcap_{s>0}\bar{\Lambda}_s=\{-\frac{1}{4}\}$, where
  $\bar{\Lambda}_s$ denotes the closure of $\Lambda_s$.
\end{itemize}
\end{theorem}

It follows that the open-arc TE and TM scattering problems~\eqref{b_conds}
can be solved by means of the second-kind integral equations
\begin{equation}\label{NSGoodDirichlet}
\tilde{N}\tilde{S}[\tilde{\varphi}]=\tilde{N}[\tilde{f}]\quad\mbox{and}
\end{equation}
\begin{equation}\label{NSGoodNeumann}
\tilde{N}\tilde{S}[\tilde{\psi}]=\tilde{g},
\end{equation}
respectively. The smooth and periodic solutions $\tilde{\varphi}$ and
$\tilde{\psi}$ of these equations are related to the singular
solutions $\mu$ and $\nu$ of equations~\eqref{bad} via
$\mu\left(\mathbf{r}(\cos\theta)\right)=\tilde{\varphi}(\theta)/\sin
\theta$ and $\nu\left(\mathbf{r}(\cos\theta)\right)=\sin\theta \left(
  \tilde{S}[\tilde{\psi}](\theta) \right).$ We thus see that, through
introduction of the weight $\omega$ and use of spaces of even and
$2\pi$-periodic functions, a picture emerges for the open-surface case
that resembles closely the one found for closed-surface
configurations: the generalized Calder\'on relation~\eqref{NSfact} is
analogous to the Calder\'on formula~\eqref{CalderonClosed}, and
mapping properties in terms of the complete range of Sobolev spaces
are recovered for $\tilde{S}$ and $\tilde{N}$, in a close analogy to
the framework presented in~\cite{Kress,Nedelec} for the operators
$\mathbf{N}_c$ and $\mathbf{S}_c$.

\subsection{Eigenvalue Distributions\label{D}}\label{sec_eigenvalues}
In order to gain additional insights on the character of the various
$k\ne 0$ open-arc operators under consideration (namely, $\tilde{S}$,
$\tilde{N}$, $\tilde{N}\tilde{S}$ as well as a generalization to
non-zero frequencies of the operator introduced in
reference~\cite{AtkinsonSloan}), we consider their corresponding
eigenvalue distributions. In Figure~\ref{EigsFigSpirale} we thus
display the eigenvalues associated with these operators, for the
spiral-shaped arc displayed in Figure~\ref{spiraleFigure} and
described in Section~\ref{sec_results}, as they were produced by means
of the quadrature rules presented in Section~\ref{sec_quadrature} and
subsequent evaluation of matrix eigenvalues. The frequency was chosen
to ensure a size to wavelength ratio $\frac{L}{\lambda}=100$, where
$L$ is the length of the arc and $\lambda$ the wavelength of the
incident wave. As expected, the eigenvalues of $\tilde{S}$ tend slowly
to zero, the eigenvalues of $\tilde{N}$ are large, while the
eigenvalues of $\tilde{N}\tilde{S}$ are bounded away from zero and
infinity, and they accumulate at $-\frac{1}{4}$.

The $k>0$ generalization of the equation~\cite{AtkinsonSloan}, whose
operator eigenvalues are displayed in the center-right portion of
Figure~\ref{EigsFigSpirale}, is obtained from right-multiplication of
the single layer operator in equation~\eqref{good} by the inverse of
the flat-arc zero-frequency single-layer operator $\tilde{S}_0^\tau$
(defined in equation~\eqref{s_tau_def} below); the resulting equation
is given by
\begin{equation}\label{AtkinsonEq}
\tilde{S}(\tilde{S}_0^\tau)^{-1}[ \tilde{\varphi}] = \tilde{f}.
\end{equation}
This equation, which can be re-expressed in the form $\left( I +(
  \tilde{S} -
  \tilde{S}_{0}^\tau)(\tilde{S}_0^\tau)^{-1}\right)[\tilde{\varphi}] =
\tilde{f}$, is a second-kind Fredholm integral equation: the operator
$\left( \tilde{S} - \tilde{S}_{0}^\tau\right)(\tilde{S}_0^\tau)^{-1}:
H^s_e(2\pi) \rightarrow H^s_e(2\pi)$ is compact.  Unfortunately, the
spectrum of the operator in equation~\eqref{AtkinsonEq} is highly
unfavorable at high frequencies, as illustrated in the center-right
image in Figure \ref{EigsFigSpirale}. Such poor spectral distributions
translate into dramatic increases, demonstrated in
Table~\ref{SloanTableSpirale}, in the number of iterations required to
solve~\eqref{AtkinsonEq} by means of Krylov subspace solvers as the
frequency grows. In fact, a direct comparison with
Table~\ref{SpiraleDirichlet} shows that the second-kind integral
equation~\eqref{AtkinsonEq} may require many more iterations at
non-zero frequencies than the original first-kind
equation~\eqref{tildeEqs}.

\section{Theoretical Considerations}\label{sec_theory}
This section presents a proof of Theorem 1. The argument proceeds by
consideration of the flat-arc zero-frequency case
(Sections~\ref{Zero-Freq} through \ref{J0sec}) and subsequent
extension to the general case (Sections \ref{Ssec},\ref{Nsec} and
\ref{NSSec}). A more detailed discussion, including full mathematical
technicalities, can be found in~\cite{LintnerBruno}.

\subsection{Flat arc at zero-frequency}\label{Zero-Freq}
In the particular case where the frequency vanishes ($k=0$) and the curve
under consideration is the flat strip ($\Gamma = [-1,1]$), the operator
$\tilde{S}$ reduces to Symm's operator
\begin{equation}\label{S0_coord}
\tilde{S}_0 [\tilde{\varphi}] (\theta ) = -\frac{1}{2\pi}\int_0^\pi
\ln|\cos\theta-\cos\theta'| \tilde{\varphi}(\theta' ) d\theta',
\end{equation}
whose well-known diagonal property~\cite{YanSloan,BrunoHaslam} in the
cosine basis $e_n(\theta) = \cos n \theta$,
\begin{equation}\label{S0Diag}
\tilde{S_0}[e_n] = \lambda_n e_n, \textbf{ }\lambda_n=
\left\{ \begin{array}{cc}\frac{\ln 2}{2}& n=0\\ \frac{1}{2n}, & n\geq 1
\end{array} \right. 
\end{equation} 
plays a central role in our analysis. Note that, in particular,
equation~\eqref{S0Diag} establishes the bicontinuity of the operator
$\tilde{S}_0$ from $H^s_e(2\pi)$ into $H^{s+1}_e(2\pi)$. The corresponding
property of bicontinuity (from $H^{s+1}_e(2\pi)$ into $H^{s}_e(2\pi)$) for
the weighted flat-strip, zero-frequency hypersingular operator
\[
\tilde{N}_0[\tilde{\psi}](\theta)=\frac{1}{2\pi}\lim\limits_{z\rightarrow 0}\int_{0}^\pi
\frac{\partial^2}{\partial z^2}\left(\ln\sqrt{(\cos\theta -
  \cos\theta')^2+z^2}\right) \tilde{\psi}(\theta')\sin^2\theta' d\theta'
\]
results, in turn, from our study of the operator $\tilde{J}_0$ in
Section~\ref{J0}.

An integration by parts argument presented in reference~\cite{Monch}
gives rise to the factorization $\tilde{N}_0 = \tilde{D}_0 \tilde{S}_0
\tilde{T}_0$
where
\begin{equation}\label{D0Def}
\tilde{D}_0[\tilde{\varphi}](\theta)=
\frac{1}{\sin\theta}\frac{d\tilde{\varphi} (\theta)}{d\theta} 
\end{equation} 
and
\begin{equation}\label{T0Def}
\tilde{T}_0 [\tilde{\varphi}](\theta)= \frac{d}{d \theta}\left(\tilde{\varphi}(\theta)
\sin\theta\right).
\end{equation}
As a result, for the flat-arc at zero-frequency case, the composite operator
$\tilde{N}\tilde{S}$, is given by
\begin{equation}\label{J0def}
\begin{split}
\tilde{J}_0=\tilde{N}_0\tilde{S}_0=\tilde{D}_0\tilde{S}_0\tilde{T}_0\tilde{S}_0.\end{split}
\end{equation}
\subsection{The operator $\tilde{J}_0$}\label{J0}

It is easy to evaluate the action of the operator $\tilde{J}_0$ on the
cosine basis: in view of equation~\eqref{J0def}, using~\eqref{S0Diag} and
expressing $\tilde{T}_0[e_n]$ as a linear combination of cosines (which
results from simple trigonometric manipulations) we obtain the relation
\begin{equation}\label{J0en}
\tilde{J}_0[e_n](\theta) = \left\{
\begin{array}{ll}-\frac{\ln2}{4}& n=0\\ 
-\cos\theta\frac{\sin n\theta}{4n\sin\theta}-\frac{\cos n\theta}{4}, & n\geq
1.
\end{array} \right.
\end{equation} 
Clearly then,
\begin{equation}\label{Jfact}
\tilde{J}_0[\tilde{\varphi}](\theta)=-\frac{\tilde{\varphi}(\theta)}{4} -
\frac{\cos\theta}{4}\tilde{C}[\tilde{\varphi}](\theta)
+\frac{1-\ln 2}{4\pi}\int_0^\pi\tilde{\varphi}(u)
du,
\end{equation}
where the operator $\tilde{C}$ is defined by
\begin{equation}\label{Hen}
\tilde{C}[e_n](\theta) = \left\{
\begin{array}{ll}0& n=0\\ 
\frac{\sin n\theta}{n\sin\theta}, & n\geq
1.
\end{array} \right.
\end{equation}

The right hand side of equation~\eqref{Jfact} may appear to bear a direct
connection with the classical closed-curve Calder\'on
formula~\eqref{CalderonClosed}. Yet, unlike the operator $\mathbf{K}_c$ in
equation~\eqref{CalderonClosed}, \emph{the operator $\tilde{C}$ is not
compact}--see Section~\ref{J0sec}. It is easy, however, to verify that
$\tilde{C}$ can also be expressed in the form 
 \begin{equation}\label{Cfact}
 \tilde{C}[\tilde{\varphi}](\theta)=\frac{\theta(\pi-\theta)}{\pi\sin\theta}\left[ \frac{1}{\theta}\int_0^\theta \tilde{\varphi}(u) du - \frac{1}{\pi-\theta}\int_\theta^\pi \tilde{\varphi}(u) du \right],
 \end{equation}
and is therefore closely related to the C\'esaro operator
$\mathcal{C}[\eta](x)=\frac{1}{x}\int_0^x \eta(u)du.$
Since $\mathcal{C}$ is bounded (but not compact) from $L^2[0,b]$ into
$L^2[0,b]$~\cite{BrownHalmosShields} (where $L^2[0,b]$ denotes the space of
square-integrable functions in the interval $[0,b]$), it follows that
$\tilde{C}: \quad H^s_e(2\pi) \rightarrow H^s_e(2\pi)$ is a continuous
operator for all $s>0$, and, therefore, $\tilde{J}_0: H^s_e(2\pi)\rightarrow
H^s_e(2\pi)$ is also a continuous operator.

Composing the operator
$\tilde{I}_0=-4\tilde{S}_0^{-1}\tilde{C}\tilde{S}_0\tilde{T}_0$ with
$\tilde{J}_0$ from the right on one hand, and from the left, on the other,
and applying the two resulting composite operators to the basis $e_n$, shows
that $\tilde{I}_0$ is the inverse of $\tilde{J}_0$, and that $\tilde{I}_0$
is a continuous operator from $H^s_e(2\pi)$ into $H^s_e(2\pi)$. It follows
that $\tilde{J}_0$ is a bicontinuous operator.  As indicated in
Section~\ref{Zero-Freq}, finally, the bicontinuity of $\tilde{N}_0$ from
$H^{s+1}_e(2\pi)$ into $H^s_e(2\pi)$ follows directly from
equation~\eqref{J0def} and the corresponding bicontinuity properties of
$\tilde{S}_0$ and $\tilde{J}_0$.

\subsection{Eigenvalues of $\tilde{J}_0$}\label{J0sec}
Re-expressing~\eqref{J0en} in the form
\begin{equation}\label{Jc0}
  \tilde J_0[e_n](\theta)=\left\{\begin{array}{ll}-\frac{\ln 2}{4}&n=0\\-\frac{\sin (n+1)\theta}{4n\sin\theta}+\frac{\cos
        n\theta}{4n}-\frac{\cos n \theta}{4},& n>0,
\end{array}\right.
\end{equation}
and making use the well-known expansion \cite{Bateman,MasonHandscomb}
\begin{equation}\label{Un}
\frac{\sin(n+1)\theta}{\sin \theta}=\left\{\begin{array}{ll} \sum\limits_{k=0}^p (2-\delta_{0k})\cos2k\theta, & n=2p \\ 2\sum\limits_{k=0}^p \cos(2k+1)\theta, & n=2p+1, \end{array}\right. 
\end{equation}
we obtain
\begin{equation}\label{J0exp}
\begin{split}
\tilde J_0[ e_n] = 
\left\{ \begin{array}{ll} \lambda_n e_n -
  \frac{1}{2n}\sum\limits_{k=0}^{p-1}(1-\frac{\delta_{0k}}{2})e_{2k},& n =
  2p,\; p\geq0\\ \lambda_n e_n - \frac{1}{2n}\sum\limits_{k=0}^{p-1}e_{2k+1},&
  n = 2p+1,\; p\geq0\\ \end{array}\right.,
  \end{split}
\end{equation}
where the diagonal elements $\lambda_n$ are given by
\begin{equation}\label{lambdaDiag}
\lambda_n=\left\{\begin{array}{lr}-\frac{\ln2}4,&n=0\\-\frac{1}{4}-\frac{1}{4n},&\quad n > 0.\end{array}\right.
\end{equation}
In view of~\eqref{Jc0}, the operator $\tilde{J}_0$ can be viewed as an
infinite upper-triangular matrix. The diagonal terms $\lambda_n$ of
this matrix are eigenvalues of $\tilde{J}_0$; the corresponding
eigenvectors $v_n$, in turn, can be expressed in terms of a finite
linear combination of the first $n$ cosine basis functions:
$v_n=\sum_{k=0}^n c_k^n e_k$.  This establishes that the set
$\Lambda_\infty$ defined in Theorem~\ref{NStheorem} is contained in
the spectrum $\sigma_s$ for all $s>0$.

As is known, the diagonal elements of an infinite upper-triangular matrix
are a subset of all the eigenvalues of the corresponding operator; for
instance, the point spectrum of the upper-triangular bounded operator
$C^*[a](n) = \sum\limits_{k=n}^\infty \frac{a_k}{k+1},$
(the adjoint of the discrete Ces\`aro operator $C$) is
given~\cite[Th. 2]{BrownHalmosShields} by the entire disc $|\lambda-1|<1$. A
similar situation arises for the operator $\tilde J_0$: searching for
sequences $\left\{ f_n, n\geq 0 \right\}$ such that $f=\sum_n f_n e_n$
satisfies $J_0[f]=\lambda f$ leads to the relations
\begin{equation}\label{frec}
  (-\frac{1}{4}-\frac{1}{4n})f_{n} -\frac{1}{2}
  \sum\limits_{k=1}^\infty \frac{f_{n+2k}}{n+2k}=\lambda f_{n}\quad,\quad   n \geq 1,
\end{equation}
and
\begin{equation}\label{frec_0}
  (-\frac{\ln 2}{4})f_{0} -\frac{1}{4}
  \sum\limits_{k=1}^\infty \frac{f_{2k}}{2k}=\lambda f_{0}\quad,\quad
  n = 0.
\end{equation}
The spectrum of $\tilde{J}_0$ equals the set of all values $\lambda$ for
which the relations~\eqref{frec} and~\eqref{frec_0} admit solutions
satisfying $\sum |f_n n^s|^2<\infty$. This set is the union of the discrete
set $\Lambda_\infty$ (see point~(iii) in Theorem~\ref{th_1}), and the open
set $\Lambda_s$ defined by
\begin{equation}\label{lambda_s} \Lambda_s=\left\{
\lambda=(\lambda_x,\lambda_y),\;4s+2 <
\frac{-\left(\lambda_x+\frac{1}{4}\right)}{(\lambda_x+\frac{1}{4})^2+\lambda_y^2}\right\}.
\end{equation}
It follows in particular that the operator $\tilde{C}$ is
not compact from $H^s_e(2\pi)$ into $H^s_e(2\pi)$: if it were, in view of
equation~\eqref{Jfact}, the spectrum of $\tilde{J}_0$ would be discrete, in
contradiction with~\eqref{lambda_s}

We thus see that, even though $\tilde{C}$ is not compact (so that, in
particular, the decomposition~\eqref{Jfact} does not present
$\tilde{J}_0$ as the sum of a multiple of the identity and a compact
operator), it follows from~\eqref{lambdaDiag} and~\eqref{lambda_s}
that the eigenvalues of $\tilde{J}_0$ are all bounded away from zero
and infinity, and that they are tightly clustered around
$-\frac{1}{4}$---in close analogy to the eigenvalue distribution
associated with the closed-surface Calder\'on
operator~\eqref{CalderonClosed}. 

\begin{table}
\caption{Scattering by a spiral-shaped arc of size $\frac{L}{\lambda} =
 800$: far-field errors.
 \label{spectralTable}}
\centering
\begin{tabular}[c]{  c  c c  c  c }
\hline 
$N$ & TE case, $\tilde{S}$ & TE case, $\tilde{N}\tilde{S}$ & TM case,
$\tilde{N}$ & TM case, $\tilde{N}\tilde{S}$ \\ 
\hline
$3000$ & $5.4\times 10^{-6}$ & $8.2\times10^{-6}$ & $1.5\times 10^{-3}$&$2.8\times 10^{-4}$ \\
$3100$ & $4.5\times 10^{-8}$ & $4.9\times 10^{-8}$ & $1.0\times 10^{-5}$  & $5.5\times 10^{-6}$ \\
$3350$ & $8.4\times 10^{-12}$ & $8.5\times10^{-12}$ & $3.7\times 10^{-10}$ & $3.7\times 10^{-11}$ \\
\hline
\end{tabular} 
\end{table}

\subsection{Properties of the operator $\tilde{S}$}\label{Ssec}
The study of $\tilde{S}$ in the general case (possibly curved arc, $k\ne 0$)
hinges on the decomposition
\begin{equation}\label{Gkernel}
G_k(\mathbf{r}(t),\mathbf{r}(t')) = A_1(k,t,t')\ln|t-t'| + A_2(k,t,t'),
\end{equation}
of the kernel $G_k$, which results from the
expression~\cite[p. 64]{ColtonKress2}
\begin{equation}\label{H0dec}
H_0^1(z)=\frac{2i}{\pi}J_0(z)\ln(z)+G(z),
\end{equation}
where $J_0(z)$ is the Bessel function and $G(z)$ is analytic. Clearly, the functions $A_1$
and $A_2$ are smooth functions of $t$ and $t'$, and the singular
behavior of $G_k$ thus resides entirely in a logarithmic term. In
particular, the continuity of the operator $\tilde{S}$ follows
directly from the corresponding continuity of the operator
$\tilde{S}_0$.

To establish the invertibility of the operator $\tilde{S}$ and the
continuity of its inverse $\tilde{S}^{-1}$, on the other hand, we
consider the decomposition
\begin{equation}\label{s_tau}
\tilde{S} =\tilde{S_0^\tau}\left( I +
\left(\tilde{S_0^\tau}\right)^{-1}(\tilde{S}-\tilde{S_0^\tau})\right),
\end{equation}
where the operator $\tilde{S}_0^\tau$, defined by 
\begin{equation}\label{s_tau_def}
\tilde{S}_0^\tau[\gamma]=\tilde{S}_0\tilde{Z}_0[\gamma ],
\end{equation}
 differs from $\tilde{S}_0$ only in the additional multiplicative term
$\tilde{Z}_0[\gamma](\theta)=\gamma(\theta)\tau(\cos\theta).$
Since the $H^s_e\to H^s_e$ operator $(\tilde{S}_0^{\tau})^{-1}(\tilde{S} -
\tilde{S}_0^\tau)$ is compact (in view of classical embedding results and
the increased regularity of the kernel of the operator
$(\tilde{S}-\tilde{S}_0^{\tau})$) and since the operator $\tilde{S}$ is
injective (as it follows from results established
in~\cite{StephanWendland,StephanWendland2,Stephan} for the
operator~\eqref{Sdef}), the bicontinuity of $\tilde{S}$ follows from the
invertibility of $\tilde{S}_0$ (and, thus, of $\tilde{S}_0^\tau$) together
with an application of the Fredholm theory to the term in parenthesis in
equation~\eqref{s_tau}.  


\subsection{Properties of the operator $\tilde{N}$}\label{Nsec}
As is known, {\em for closed curves}, the operator $\mathbf{N}_c$ (the
normal derivative of the double-layer potential), can be expressed in terms
of a composition of a single layer potential and derivatives tangential to
the curve $\Gamma_c$~\cite[p. 57]{ColtonKress1}. While for open surfaces the
conversion of normal derivatives into tangential derivatives for the
operator $\mathbf{N}$ gives rise to highly singular boundary terms (which
arise from the integration-by-parts calculation that is part of the
conversion process), the presence of the boundary-vanishing weight $\omega$
in our $N_\omega$ operator eliminates all singular terms, and we obtain
$N_\omega[\varphi](t)=N_\omega^g[\varphi](t)+N_\omega^{pv}[\varphi](t)$ %
where
\begin{equation}\label{NomegaG}
N_\omega^g[\varphi](t)=k^2\int_{-1}^1 G_k(\mathbf{r}(t),\mathbf{r}(t'))\;\varphi(t')\; \tau(t')\;\sqrt{1-t'^2}\;\mathbf{n}_t \cdot \mathbf{n}_{t'}\;dt',\end{equation}
and where
\begin{equation}\label{NomegaPV}
  N_\omega^{pv}[\varphi](t)=\frac{1}{\tau(t)}\frac{d}{dt}\int_{-1}^1
    G_k(\mathbf{r}(t),\mathbf{r}(t'))\;\frac{d}{dt'} \left(
      \varphi(t')\;\sqrt{1-t'^2}\right)dt',
\end{equation}
see~\cite{Monch}. Using the changes of variables $t=\cos\theta$ and
$t'=\cos\theta'$ in equations~\eqref{NomegaG} and~\eqref{NomegaPV},
 it follows that
\begin{equation}\label{Nfact}
  \tilde{N}[\tilde\varphi](\theta)=\tilde{N}^{g}[\tilde\varphi](\theta)+\tilde{N}^{pv}[\tilde\varphi](\theta),
\end{equation}
where
\begin{equation}\label{Ng}
\begin{split}
\tilde{N}^g[\tilde\varphi](\theta)=k^2\int_{0}^\pi
G_k(\mathbf{r}(\cos\theta),\mathbf{r}(\cos\theta'))
\;\tilde\varphi(\theta')\; \tau(\cos\theta')\;\\\sin^2\theta'
\mathbf{n}_\theta \cdot \mathbf{n}_{\theta'}\;d\theta',
\end{split}
\end{equation}
and where, defining $\tilde
T_0^\tau[\tilde\varphi](\theta)=\frac{1}{\tau(\cos\theta)}T_0[\tilde
\varphi](\theta)$,
\begin{equation}\label{Npv}
\tilde{N}^{pv}[\tilde\varphi](\theta) =
\frac{1}{\tau(\cos\theta)}\left(\tilde D_0 \tilde{S} \tilde
T_0^\tau\right)[\tilde\varphi](\theta).
\end{equation}
Note for future use that equation~\eqref{Npv} can be re-expressed in the
form
\begin{equation}\label{pv_0}
{\tilde N}^{pv}[\tilde \varphi](\theta) =\tilde{N}^\tau_0[\tilde
\varphi](\theta) + \frac{1}{\tau(\cos\theta)} \tilde D_0
(\tilde{S}-\tilde{S}_0^\tau) \tilde T_0^\tau [\tilde \varphi](\theta).
\end{equation}
where $\tilde{N}^\tau_0[\gamma] = \tilde{Z}_0^{-1}\tilde{N}_0[\gamma].$

In view of equation~\eqref{Nfact}, the operator ${\tilde N}$ equals the sum
of $\tilde{N}^g$ (which, like $\tilde{S}$, maps $H^s_e(2\pi)$ into
$H^{s+1}_e(2\pi)$) and $\tilde{N}^{pv}$. But, from~\eqref{pv_0}, we see that
$\tilde{N}^{pv}$ can be expressed as the sum of $\tilde{N}_0^{\tau}$ and a
compact perturbation. Since, like $\tilde{N}_0$, the operator
$\tilde{N}_0^{\tau}$ maps $H^{s+1}_e(2\pi)$ into $H^s_e(2\pi)$, it follows
that $\tilde{N}^{pv}$ is a bounded operator from $H^{s+1}_e(2\pi)$ into
$H^s_e(2\pi)$, and, thus, $\tilde{N}: \quad H^{s+1}_e(2\pi) \rightarrow
H^s_e(2\pi)$ is a bounded operator as well.

The bicontinuity properties of $\tilde{N}$ can now be established by
an argument similar to the one applied to $\tilde{S}$ in the previous
section using, this time, the identity
$\tilde{N}=\tilde{N^\tau_0}\left( I
  +\left(\tilde{N}^\tau_0\right)^{-1}(\tilde{N}-\tilde{N}^\tau_0)\right).$

%

\begin{figure}[h]
\begin{center}
\includegraphics[scale=0.21]{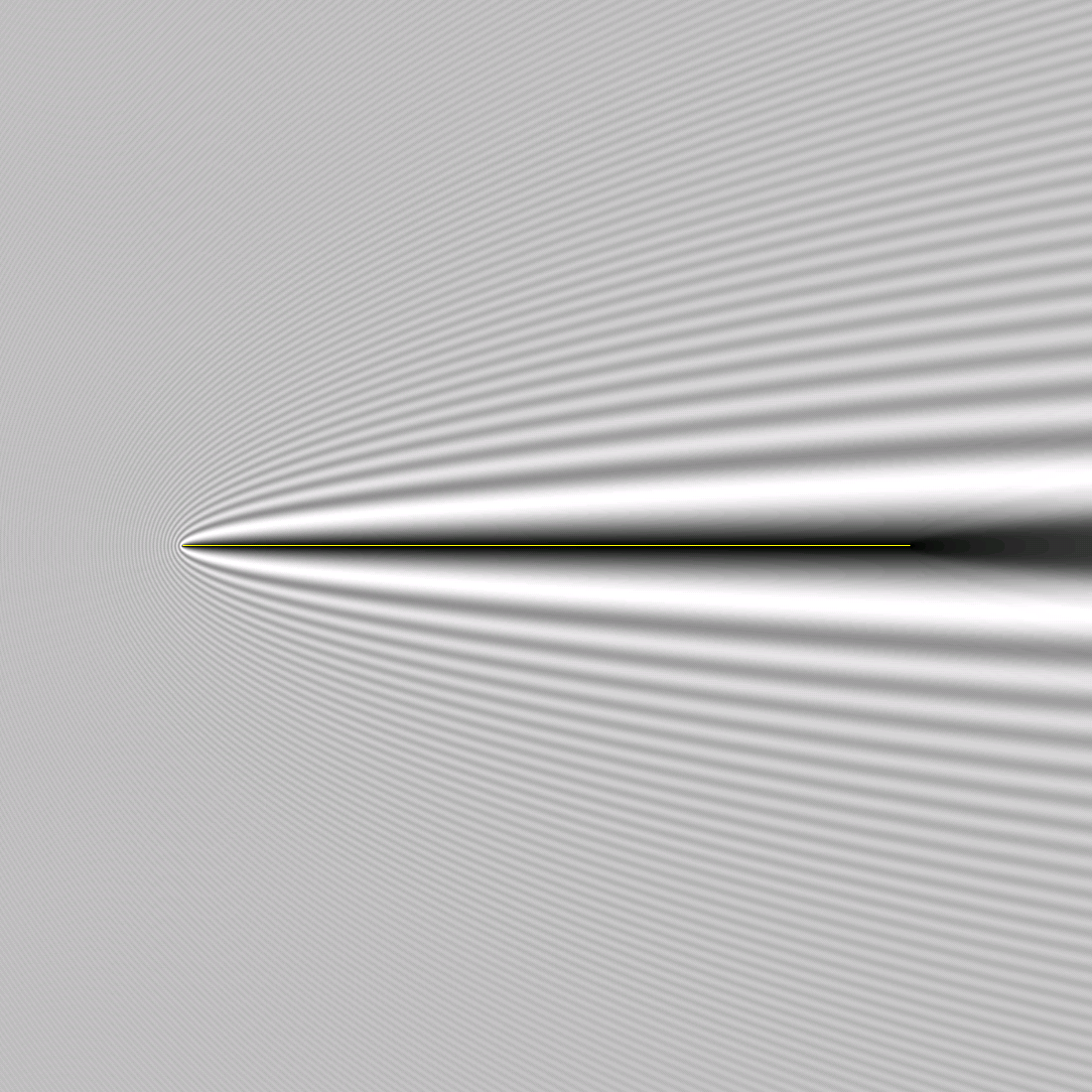}\hspace{0.5cm}
\includegraphics[scale=0.375]{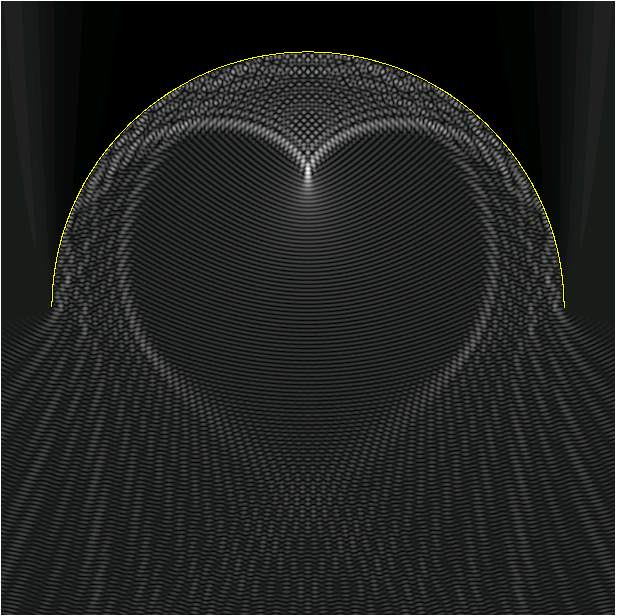}
\caption{Sample diffraction patterns. Left: TE scattering (total field) by
  an infinitely thin strip of size $L=200\lambda$ for horizontal
  left-to-right incidence. (The TM total field at this incidence equals the
  incident field.)  Right: TM solution for a half-circle of size
  $L=200\lambda$ under normal incidence from below.  Note the trailing
  shadow in the wake of the strip, and the caustics inside the circular
  reflector.
\label{Strip_Figure}} 
\end{center}
\end{figure}

\subsection{Generalized Calder\'on formula}\label{NSSec}
The generalized Calder\'on formula is obtained by expressing the combination
$\tilde{N}\tilde{S}$ in the form
\begin{equation}\label{CalderonFinal}
\tilde{N}\tilde{S}=\tilde{N}_0^\tau\tilde{S}_0^\tau + \tilde K
=\tilde{J}_0^\tau + \tilde K
\end{equation}
where the operator $\tilde K$, given by $\tilde K=\tilde{N}(\tilde{S} -
\tilde{S}_0^\tau) + (\tilde{N}- \tilde{N}_0^\tau)\tilde{S}_0^\tau,$
is compact in view of the classical Sobolev embedding theorems and the
(increased) smoothness properties of the kernels of the operators
$(\tilde{S}-\tilde{S}_0^\tau)$ and $(\tilde{N}- \tilde{N}_0^\tau)$
respectively.  In view of the relation $\tilde{J}_0^\tau =
\tilde{Z}_0^{-1}\tilde{J}_0\tilde{Z}_0,$
$\tilde{J}_0$ and $\tilde{J}_0^\tau$ have the same spectrum:
$(\lambda, v)$ is an eigenvalue-eigenvector pair for $\tilde{J}_0$ if
and only if $(\lambda, \tilde{Z}_0^{-1}[v])$ is an
eigenvalue-eigenvector pair for $\tilde{J}_0^\tau$. Clearly,
$\tilde{J}_0^\tau$ has the same mapping properties and regularity as
$\tilde{J}_0$ ($\tilde{J}_0^\tau$ is bicontinuous from $H^s_e(2\pi)$
into $H^s_e(2\pi)$), and therefore, the
decomposition~\eqref{CalderonFinal} shows that $\tilde{N}\tilde{S}$
equals the sum of an invertible bicontinuous operator and a compact
operator; that is, $\tilde{N}\tilde{S}$ is a second-kind Fredholm
operator.

\section{High-Order Numerical Methods\label{sec_quadrature}}
In this section we present spectral quadrature rules for the operators
$\tilde{S}$, $\tilde{N}$ and $\tilde{N}\tilde{S}$ which give rise to an
efficient and accurate solver for the general open arc diffraction
problem~\eqref{b_conds}.
\subsection{Spectral discretization for $\tilde{S}$\label{spec_S}}
Use of the nodes $\left\{\theta_n=\frac{\pi(2n+1)}{2N}\right\}$,
${n=0,\cdots, N-1}$, gives rise \cite[eq. (5.8.7),(5.8.8)]{NumericalRecipes}
to a spectrally convergent cosine
representation for smooth, $\pi$-periodic and even functions
$\tilde{\varphi}$: \begin{equation}\label{chebexp}
\begin{split}
\tilde{\varphi}(\theta)=\sum_{n=0}^{N-1} a_n \cos( n\theta
),\quad\mbox{where}\quad
a_n=\frac{(2-\delta_{0n})}{N}\sum\limits_{j=0}^{N-1} \tilde{\varphi}(
\theta_j) \cos (n \theta_j).
\end{split}
\end{equation}
Thus, applying equation~\eqref{S0Diag} to each term of
expansion~\eqref{chebexp}, we obtain the well-known \emph{spectral}
quadrature rule for the logarithmic kernel
\begin{equation}\label{logrule}
\int_0^\pi \ln|\cos\theta-\cos\theta'|\tilde{\varphi}(\theta')d\theta'\sim \frac{\pi}{N}\sum\limits_{j=0}^{N-1} \tilde{\varphi}(\theta_j)R^{(N)}_j(\theta),\end{equation}
where \begin{equation}\label{Rn}
R^{(N)}_j(\theta)=-2\sum\limits_{m=0}^{N-1}(2-\delta_m)\lambda_m\cos(m\theta_j)\cos(m\theta).
\end{equation}
Following~\cite{ColtonKress2,Martensen,Kussmaul} we then devise a high-order
integration rule for the operator $\tilde{S}$, noting first
from~\eqref{Gk_def} and~\eqref{Gkernel}
that
\begin{equation}\label{Gk2}
\begin{split}
G_k(\mathbf{r}(\theta),\mathbf{r}(\theta'))=
A_1(k,\cos\theta,\cos\theta')\ln|\cos\theta-\cos\theta'| + A_2(k,\cos\theta,\cos\theta'),
\end{split}
\end{equation}
where, letting $R = |\mathbf{r}(\cos\theta)-\mathbf{r}(\cos\theta')|$
we have
\begin{equation}
A_1(k,\cos\theta,\cos\theta')=-\frac{1}{2\pi}J_0(kR),
\end{equation} and
\begin{equation}
A_2(k,\cos\theta,\cos\theta')=\frac{i}{4}H^{1}_0(kR)+\frac{1}{2\pi}J_0(kR)\ln|\cos\theta-\cos\theta'|.
\end{equation}
In view of~\eqref{H0dec} and the smoothness of the ratio
$\frac{R}{|\cos\theta-\cos\theta'|}$, the functions $A_1$ and $A_2$ are
even, smooth (analytic for analytic arcs) and $2\pi$-periodic functions of
$\theta$ and $\theta'$---and, thus, in view of~\eqref{logrule}, the
expression
\begin{equation}
\begin{split}
\int_0^\pi
\tilde{\varphi}(\theta')A_1(k,\cos\theta,\cos\theta')\ln|\cos\theta-\cos\theta
'|\tau(\cos\theta')d\theta' \\ \sim
\frac{\pi}{N}\sum_{j=0}^{N-1}\tilde\varphi(\theta_j)\tau(\cos\theta_j)A_1(k,\cos\theta,\cos\theta_j)R^{(N)}_j(\theta)
\end{split}
\end{equation}
provides a spectrally accurate quadrature rule. By making use of trapezoidal
integration for the second term in the right-hand side of~\eqref{Gk2} we therefore
obtain the spectrally accurate quadrature approximation of the operator
$\tilde{S}$\,:
\begin{equation}
\begin{split}
\tilde{S}[\varphi](\theta)\sim  \frac{\pi}{N}\sum_{j=0}^{N-1}\tilde\varphi(\theta_j)\tau(\cos\theta_j)
  \big(A_1(k,\cos\theta,\cos\theta_j)R^{(N)}_j(\theta)
  +A_2(k,\cos\theta,\cos\theta_j) \big)\label{Srule2D_2}.
\end{split}
\end{equation}

\subsection{Efficient implementation}
The right-hand side of~\eqref{Srule2D_2} can be evaluated directly for all
$\theta$ in the set of quadrature points
$\left\{\theta_n,n=0,\dots,N-1\right\}$ by means of a matrix-vector
multiplication involving the matrix $S^{(N)}$ whose elements are defined by
\begin{equation}\label{Smatrix}
\begin{split}
S^{(N)}_{nj}=\frac{\pi}{N}\tau(\cos\theta_j)\left( A_1(k,\cos\theta_n,\cos\theta_j)R^{(N)}_j(\theta_n)+A_2(k,\cos\theta_n,\cos\theta_j)\right).
\end{split}
\end{equation}
A direct evaluation of the matrix $S^{(N)}$ on the basis of~\eqref{Rn}
requires $O(N^3)$ operations; as shown in what follows, however, the matrix
$S^{(N)}$ can be produced at significantly lower computational cost. Indeed,
expressing the product of cosines in~\eqref{Rn} as a sum of cosines of added
and subtracted angles, the quantities
\begin{equation}
R^{(N)}_j(\theta_n)=-\sum_{m=0}^{N-1}(2-\delta_m)\lambda_m\big(\cos(\frac{m\pi}{N}|n-j|)+\cos(\frac{m\pi}{N}(n+j+1))\big)
\end{equation}
can be expressed in the form 
\begin{equation}
R^N_j(\theta_n)= R^{(N)}(|n-j|) + R^{(N)}(n+j+1),
\end{equation}
where the vector $R^{(N)}$ is given by 
\begin{equation}
R^{(N)}(\ell)=
-\sum_{m=0}^{N-1}(2-\delta_m)\lambda_m\cos(\frac{m\pi}{N}\ell),\textbf{
}\ell\in[0,2N-1].
\end{equation}
Our algorithm evaluates this vector efficiently by means of an FFT, and
produces as a result the matrix $S^{(N)}$ at an overall computational cost
of $O(N^2\ln N)$ operations. This fast spectrally-accurate algorithm could
be further accelerated, if necessary, by means of techniques such as those
presented in
References~\cite{BleszynskiBleszynskiJaroszewicz,BrunoKunyansky,Rokhlin}.

\begin{figure}
\begin{center}
\includegraphics[scale=0.24]{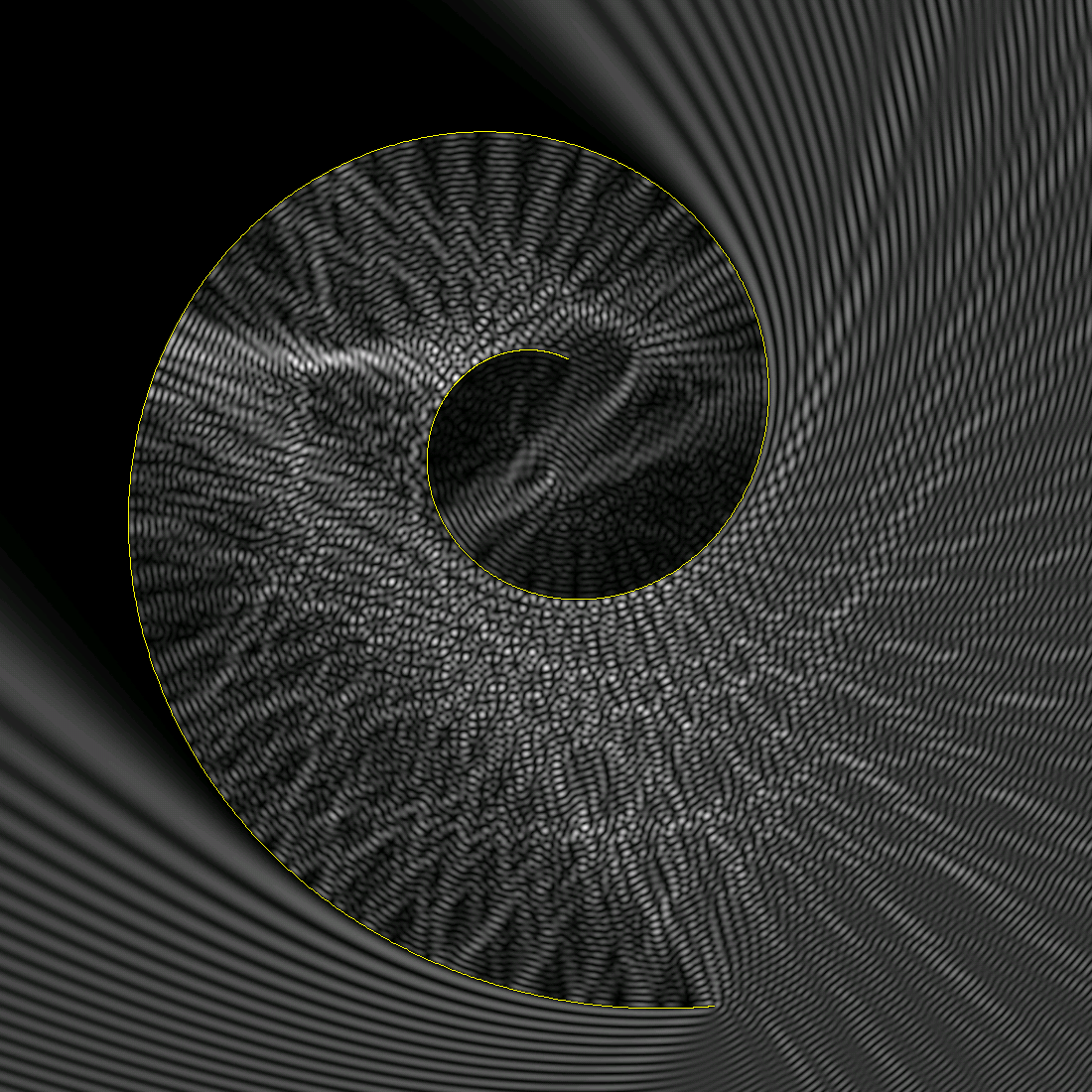}\hspace{0.5cm}
\includegraphics[scale=0.24]{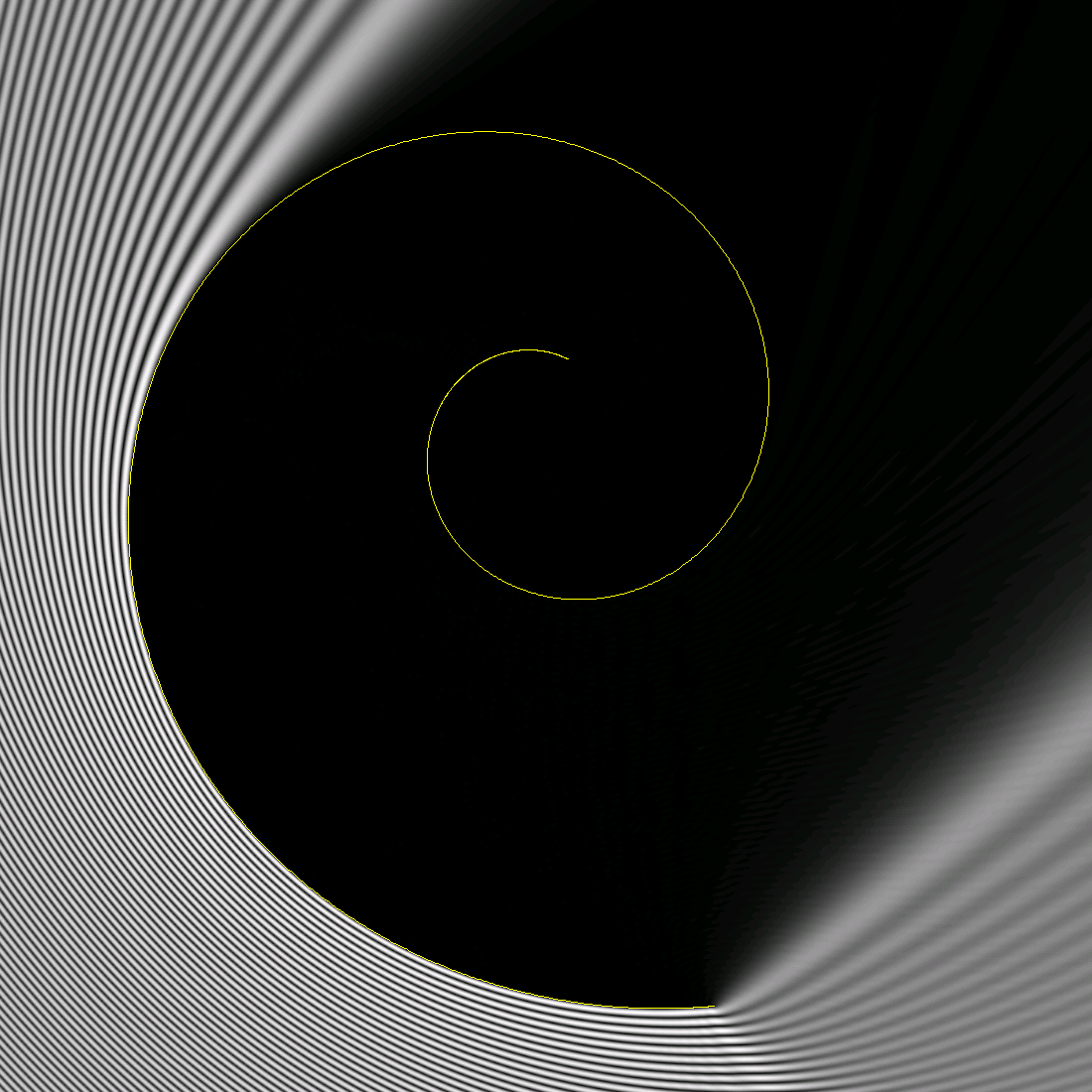}
\caption{TE-diffraction by a spiral-shaped arc of size $L=200\lambda$, for
  incidence angles of 135$^\circ$ (left) and 45$^\circ$ (right) from the
  positive $x$-axis. The left figure shows internal reflections that enable
  the field to penetrate to the center of the spiral, giving rise to an
  interesting array of caustics. The relative error $\epsilon_r$ in both
  numerical solutions is no larger than $10^{-5}$.
 \label{spiraleFigure}} 
\end{center}
\end{figure}

\subsection{Spectral Discretization for $\tilde{N}$}
In order to evaluate $\tilde{N}$ we use~\eqref{Nfact}, the first term of which
is a single-layer operator which can be evaluated by means of a rule
analogous to~\eqref{Srule2D_2} and a rapidly computable matrix $N^{g,(N)}$ (similar to $S^{(N)}$)
with elements
\begin{equation}
N^{g,(N)}_{nj}=\frac{k^2\pi}{N}\tau(\cos\theta_j)\sin^2\theta_j
(\mathbf{n}_{\theta_j}\cdot \mathbf{n}_{\theta_n})
\left(A_1(k,\cos\theta_n,\cos\theta_j)R^{(N)}_j(\theta_n)+A_2(k,\cos\theta_n,\cos\theta_j)\right).
\end{equation}
To evaluate the second term in~\eqref{Nfact}, in turn, we make use of the
decomposition~\eqref{Npv}, and we approximate the quantity
$\tilde{T}_0^\tau[\tilde{\varphi}]$ by means of term per term
differentiation of the sine expansion of the function
$\tilde{\varphi}(\theta)\sin(\theta)$ (which can itself be produced
efficiently by means of an FFT). Since $\tilde{D}_0$ is essentially the
differentiation operator in the $x$ variable,
\[\frac{1}{\sin\theta}\frac{d}{d\theta}\left(
\varphi(\cos\theta)\right)=-\frac{d}{dx}\left(\varphi(x)\right),
\] 
our solver evaluates the quantity $\tilde{D}_0[\tilde{\varphi}]$ by invoking
classical FFT-based Chebyshev differentiation
rules~\cite[p. 195]{NumericalRecipes}.

 \begin{table}
\begin{center}
\begin{tabular}[c] {  c  c c  c  c  c c  c }

\multicolumn{4}{ c  }{} & \multicolumn{2}{ c }{ TE($\tilde{S}$)} & \multicolumn{2}{ c }{TE($\tilde{N}\tilde{S}$)}  \\
  \hline
  $\frac{L}{\lambda}$ & N  & $\epsilon_r$ & Mat.   &  It. & Time &  It. & Time \\
\cline{1-8}
\cline{1-8}
50& 400 &$<10^{-5}$& $<1s$ & 24  &$<1s$ & 8 & $<1s$   \\

200& 1600 &$<10^{-5}$&$4s$ & 33 & $1s$ &8 & $2s$\\

800& 6400 &  $<10^{-5}$&$54s$ &45 & $18s$ & 8 & $ 15s$\\

\hline
\end{tabular}

 \caption{ Conditioning and computing times for the TE (Dirichlet) problem
 on the flat strip.\label{StripDirichlet}}

\end{center}
\end{table}

 \begin{table}
\begin{center}
\begin{tabular}[c] {  c  c c  c  c  c c  c }
\multicolumn{4}{ c  }{} & \multicolumn{2}{ c }{ TE($\tilde{S}$)} & \multicolumn{2}{ c }{TE($\tilde{N}\tilde{S}$)}  \\
  \hline
  $\frac{L}{\lambda}$ & N  & $\epsilon_r$ & Mat.   &  It. & Time &  It. & Time \\
\cline{1-8}
\cline{1-8}
50& 400 &$<10^{-5}$& $<1s$ & 64  &$<1s$  & 46 & $<1s$ \\

200& 1600 &$<10^{-5}$&$4s$ & 93  & $3s$  & 62 & $8s$\\

800& 6400 &$<10^{-5}$& $55s$& 136 & $58s$ & 79 & $158s$\\
\hline
\end{tabular}

 \caption{ Conditioning and computing times for the TE (Dirichlet) problem
 on the spiral-shaped arc.\label{SpiraleDirichlet}}
\end{center}
\end{table}

\begin{table}
\begin{center}
\begin{tabular}[c] {  c  c c  c  c  c c  c }
\multicolumn{4}{ c }{} & \multicolumn{2}{ c }{ TM($\tilde{N})$} &
  \multicolumn{2}{ c }{TM($\tilde{N}\tilde{S})$} \\ \hline $\frac{L}{\lambda}$ &
  N & $\epsilon_r$ & Mat.  & It. & Time & It. & Time \\ \cline{1-8}
  \cline{1-8} 50& 400 & $< 10^{-5}$& $<1s$ & 67 & $<1s$ &9 &$<1s$ \\
200& 1600 &$< 10^{-5}$ & $4s$ & 160 & $16s$ & 9 &  $1s$ \\
800& 6400 & $< 10^{-5}$ & $55s$ & 298  &$415s$ & 9 & $17s$\\
 

\hline
\end{tabular}
\caption{Conditioning and computing times for the TM (Neumann) problem on
the strip.
\label{StripNeumann}}
\end{center}
\end{table}

 \begin{table}
\begin{center}
\begin{tabular}[c] {  c  c c  c  c  c c  c }

\multicolumn{4}{ c  }{} & \multicolumn{2}{ c }{ TM($\tilde{N})$} & \multicolumn{2}{ c }{TM($\tilde{N}\tilde{S}$)}  \\
  \hline
  $\frac{L}{\lambda}$ & N  & $\epsilon_r$ & Mat.   &  It. & Time &  It. & Time \\
\cline{1-8}
\cline{1-8}
50& 400 &$<10^{-5}$& $<1s$ & 202 &$<1s$ & 48 & $<1s$\\
200& 1600 &$<10^{-5}$&$3s$ & 432 & $65s$ & 63 & $8s$  \\
800& 6400 &$<10^{-5}$& $55s$ & 849  & $1692s$ & 83 & $160s$\\

\hline
\end{tabular}

\caption{ Conditioning and computing times for the TM (Neumann) problem on
the spiral-shaped arc. \label{SpiraleNeumann}}
\end{center}
\end{table}

\begin{figure}
\begin{center}
\vspace{2mm}
\includegraphics[scale=0.32]{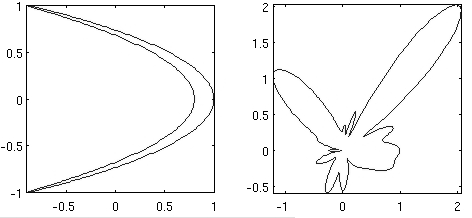}
\includegraphics[scale=0.32]{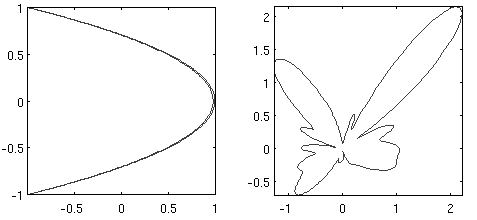}
\includegraphics[scale=0.32]{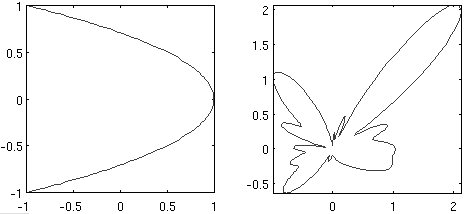}

\caption{A sequence of increasingly thin closed curves converging to the
  open parabolic scatterer $x=1-2y^2$ and corresponding far-field patterns.
  From left to right: closed-curve with $a$=0.9, corresponding far field,
  closed-curve with $a$=0.99, corresponding far field, parabolic (open) arc
  and corresponding far field. Note the convergence of the far-field
  patterns as the closed scatterers approach the open parabolic scatterer.
  \label{KiteFigure}}
\end{center}
\end{figure}

\section{Numerical Results }\label{sec_results}
The numerical results presented in what follows were obtained by means of a
C++ implementation of the quadrature rules introduced in
Section~\ref{sec_quadrature} for numerical evaluation of the operators
$\tilde{S}$ and $\tilde{N}$ (and thus, through composition,
$\tilde{N}\tilde{S}$), in conjunction with the iterative linear algebra
solver GMRES~\cite{saad}.  In all cases the errors reported were evaluated
by comparisons with highly-resolved numerical solutions. All runs were
performed in a single 2.2GHz Intel processor. 

\noindent
{\bf Spectral convergence.} 
To demonstrate the high-order character of the algorithm described in
previous sections we consider the problems of TE and TM scattering by
the exponential spiral 
$x(s)=e^s \cos(5s)$, $y(s)=e^s\sin(5s)$
of size $\frac{L}{\lambda}=800$, where $L$ and $\lambda$ denote the
perimeter of the curve and the electromagnetic wavelength, respectively.
Table~\ref{spectralTable} demonstrates the spectral (exponentially-fast)
convergence of the TE and TM numerical solutions produced by means of the
operators $\tilde{S}$, $\tilde{N}$ and $\tilde{N}\tilde{S}$ for this problem
(cf. equations~\eqref{tildeEqs} and~\eqref{NSGoodDirichlet}
and~\eqref{NSGoodNeumann}); note from Figure~\ref{spiraleFigure} the
manifold caustics and multiple reflections associated with this solution.

\noindent
{\bf Limit of Closed Curves.} 
In order to obtain an indication of the manner in which an open arc problem
can be viewed as a limit of closed-curve problems (and, in addition, to
provide an independent verification of the validity of our solvers) we
consider a test case in which the open-arc parabolic scatterer $x=1-2y^2$ is
viewed as the limit as $a \rightarrow 1$ of the family of closed curves
$x(s)=(1-a)\cos s+a\cos(2s)$, $y(s)=\sin(s).$
Using the closed-curve Nystr\"om algorithms~\cite{ColtonKress2} we evaluate
the TE fields scattered by these closed curves at $k=10$ for values of $a$
approaching $a=1$. Figure~\ref{KiteFigure} displays the $k=10$ far fields
corresponding to $a=0.9$ and $a=0.99$ side-by-side the corresponding
far-field pattern for the limiting open parabolic arc as produced by the
$\tilde{S}$-based open-arc solver. Clearly the closed-curve and open-arc
solutions are quite close to each other.  As might be expected, as $a$
approaches 1 an increasingly dense discretization is needed to maintain
accuracy in the closed-curve solution: for $a=0.9$, 256 points where needed
to reach a far field error of $10^{-4}$, while for $a=0.99$ as many 1024
points were needed to reach the same accuracy---even for the low frequency
under consideration. The corresponding open-arc solution, in contrast, was
produced with $10^{-4}$ accuracy by means of a much coarser, $64$ point
discretization.

\noindent
{\bf Solver performance.}  The TE (Dirichlet) problem can be solved by means
of either the left-hand equation in~\eqref{tildeEqs} or
equation~\eqref{NSGoodDirichlet}, which in what follows are called equations
TE$(\tilde{S})$ and TE$(\tilde{N}\tilde{S})$, respectively. The TM (Neumann)
problem, similarly, can be tackled by means of either the right-hand
equation in~\eqref{tildeEqs} or equation~\eqref{NSGoodNeumann}; we call
these equations TM$(\tilde{N})$ and TM$(\tilde{N}\tilde{S})$, respectively.
Results for TE and TM problems obtained by the various relevant equations
for two representative geometries, a strip $[-1,1]$ and the exponential
spiral mentioned above in this section, are presented in
Figures~\ref{Strip_Figure} and~\ref{spiraleFigure} and
Tables~\ref{StripDirichlet} through~\ref{SpiraleNeumann}. In the tables the
abbreviation ``It.'' denotes the number of iterations required to achieve an
$\epsilon_r$ relative maximum error in the far field (calculated as the
quotient of the maximum absolute error in the far field by the maximum
absolute value of the far-field), ``Mat.''  is the time needed to build the
$\tilde S$ matrix given by equation~\eqref{Smatrix}, as well as (when
required) the corresponding matrix for $\tilde{N}^g$ which can be
constructed in $O(N^2)$ operations from the $\tilde{S}$ matrix (see
Section~\ref{sec_quadrature}), and ``Time'' is the total time required by
the solver to find the solution once the matrix is stored.

As can be seen from these tables, the TM equation TM$(\tilde{N})$ requires
very large number of iterations as the frequency grows and, thus, the
computing times required by the low-iteration second-kind equation
TM$(\tilde{N}\tilde{S})$ are significantly lower than those required by
TM$(\tilde{N})$. The situation is reversed for the TE problem: although, the
corresponding second-kind equation TE$(\tilde{N}\tilde{S})$ requires fewer
iterations than TE$(\tilde{S})$, the total computational cost of the
second-kind equation is generally higher in this case---since the
application of the operator in TE$(\tilde{S})$ is significantly less
expensive than the application of operator in TE$(\tilde{N}\tilde{S})$.

\begin{figure}
\begin{center}
\includegraphics[scale=0.22]{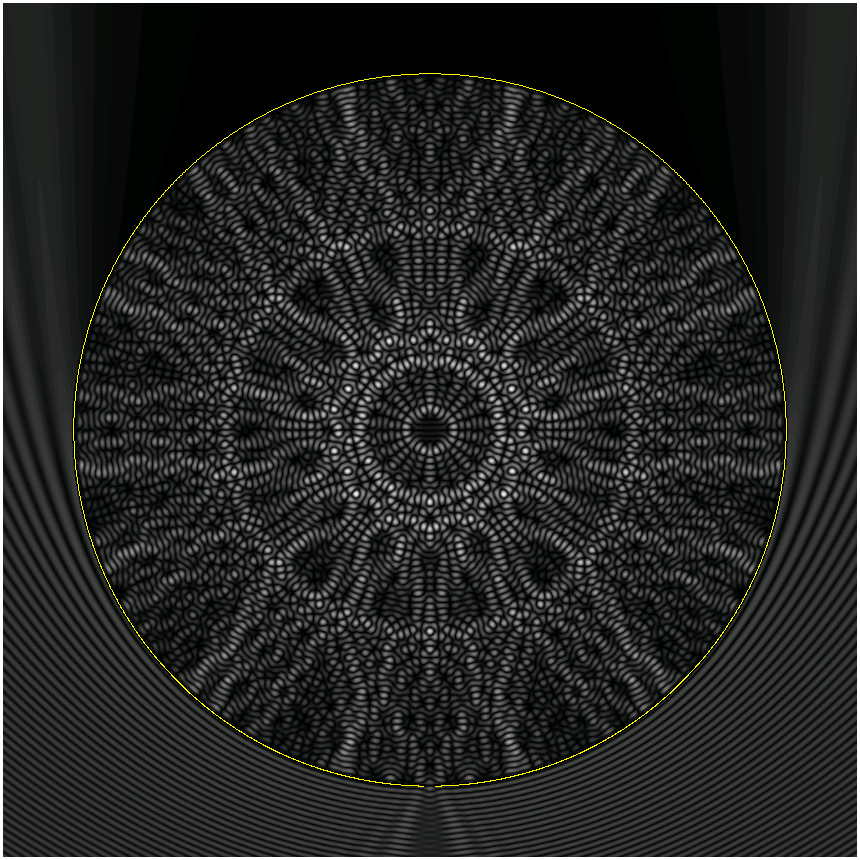}
\includegraphics[scale=0.22]{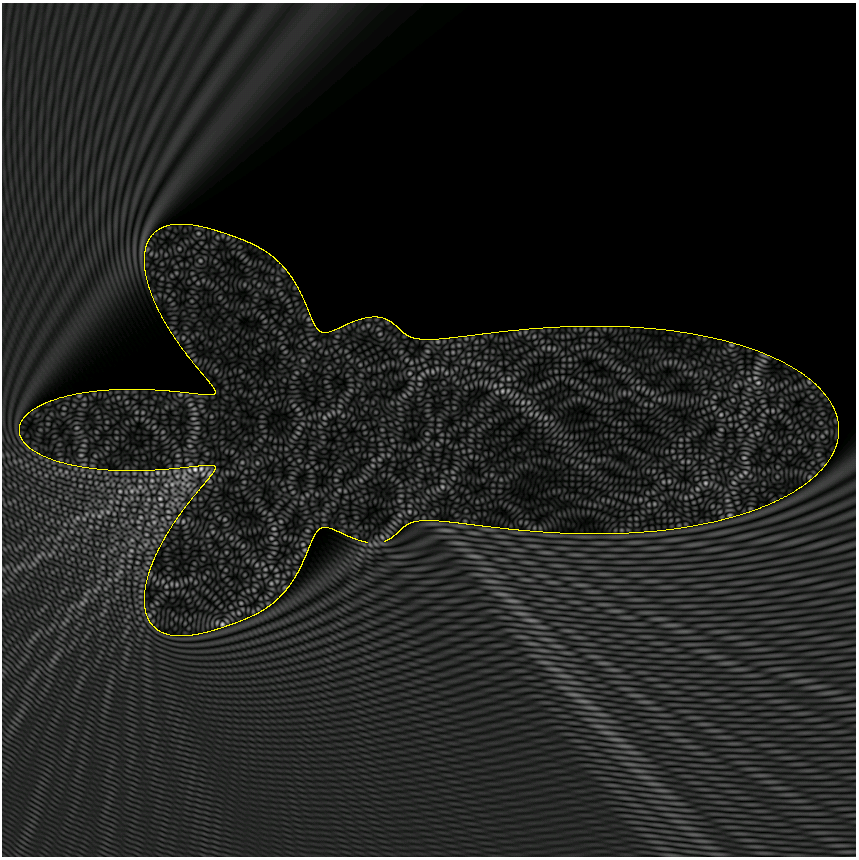}
\includegraphics[scale=0.22]{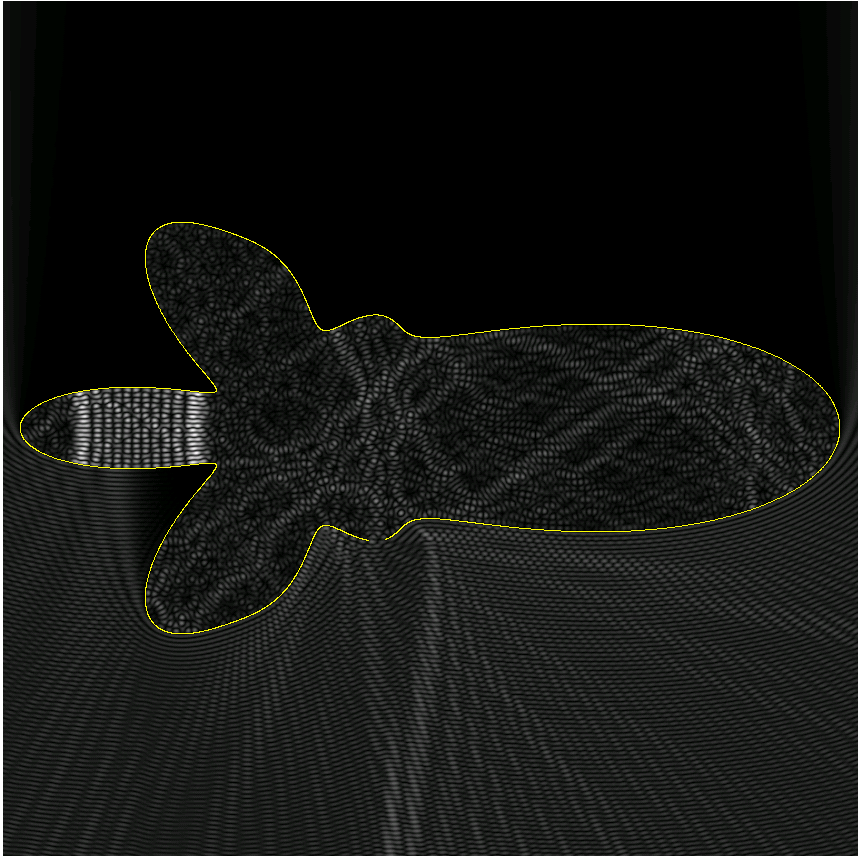}
\caption{ 
Left: TM polarization, for a circular cavity of size
  $\frac{L}{\lambda}=200$, with a bottom aperture of size equal to the
  wavelength $\lambda$.  Center: TE problem for a rocket-shaped cavity with
  perimeter $L=200\lambda$, and a bottom aperture equal to
  $6\lambda$. Right: same as in the center figure, but in the case of normal
  incidence---for which a strong resonance develops in the tail of the
  rocket.
  \label{Circ_res}} 
\end{center}
\end{figure}




\noindent {\bf Resonant Cavities.} We have found that interesting
resonant electromagnetic behavior arises from diffractive elements
constructed as almost-closed open-arcs. As can be seen in
Figure~\ref{Circ_res}, for example, circular and rocket-shaped
cavities with small openings (a few wavelengths in size), can give rise to
interesting and highly energetic field patterns within the open
cavity.  The number of iterations required for each of these
configurations is of course much larger than for simpler geometries,
such as the strip.  Yet, overall reduction in number of iterations and
computing times is observed when the equation TM$(\tilde{N}\tilde{S})$
is used in lieu of TM$(\tilde{N})$. For the TE problem, once again
TE$(\tilde{S})$ gives rise to faster overall numerics than
TE$(\tilde{N}\tilde{S})$, although the latter equation still requires
fewer GMRES iterations.

\begin{table}
\begin{center}
\begin{tabular}[c]{ c  c c  c  c  c}
\hline
$\frac{L}{\lambda}$& N  &$\epsilon_r$& It. & Mat. Time &Sol. Time \\ 
\hline
50& 400 &$1.2\times10^{-5}$& 124 & $<1s$ &$ 1s$ \\
200& 1600 &$6.3\times10^{-6}$& 293  &$3s$ & $15s$ \\
800& 6400 &$2.2\times10^{-5}$& 672 & $56s$ & $411s$ \\
\hline
\end{tabular}
\caption{Conditioning and times for the solution of the Dirichlet problem using the
  generalization of the method~\cite{AtkinsonSloan,RokhlinJiang}:
  $(\tilde{S}_{0}^\tau)^{-1}$ as a preconditioner for the spiral-shaped arc. \label{SloanTableSpirale}}
\end{center}
\end{table}

\noindent {\bf Conclusions.} We have introduced new second kind equations
and associated numerical algorithms for solution of TE and TM scattering
problems by open arcs. The new open-arc second-kind formulations are the
first ones in the literature that lead to reduced number of GMRES iterations
consistently across various geometries and frequency regimes.  The new
second-kind formulation TM$(\tilde{N}\tilde{S})$ is highly beneficial for
the TM problem, giving rise to order-of-magnitude improvements in computing
times over the hypersingular formulation TM$(\tilde{N})$. Such gains do not
occur in the TE case: although, the second-kind equation
TE$(\tilde{N}\tilde{S})$ requires fewer iterations than TE$(\tilde{S})$, the
total computational cost of the second-kind equation is generally higher in
this case---since the application of the operator in TE$(\tilde{S})$ is
significantly less expensive than the application of operator in
TE$(\tilde{N}\tilde{S})$. Generalization of these methods enabling efficient
second-kind solution of problems of scattering by open-surfaces in three
dimensions will be presented elsewhere~\cite{BrunoLintner3}.

\noindent
{\bf Acknowledgments.} The authors gratefully acknowledge support from
  AFOSR, NSF, JPL and the Betty and Gordon Moore Foundation.

\bibliographystyle{plain} \bibliography{./Screens}

\begin{thebibliography}{10}

\bibitem{AntoineBendaliDarbasMarion}
Xavier Antoine, Abderrahmane Bendali, and Marion Darbas.
\newblock Analytic preconditioners for the boundary integral solution of the
  scattering of acoustic waves by open surfaces.
\newblock {\em J. Comput. Acoust.}, 13(3):477--498, 2005.

\bibitem{AtkinsonSloan}
K.~Atkinson and I.~Sloan.
\newblock The numerical solution of first-kind logarithmic-kernel integral
  equations on smooth open arcs.
\newblock {\em Math. Comp.}, 56(193):119--139, 1991.

\bibitem{BleszynskiBleszynskiJaroszewicz}
E~Bleszynski, M~Bleszynski, and T~Jaroszewicz.
\newblock {AIM: Adaptive integral method for solving large-scale
  electromagnetic scattering and radiation problems}.
\newblock {\em {Radio Sci.}}, {31}({5}):{1225--1251}, {1996}.

\bibitem{BrownHalmosShields}
A.~Brown, P.~R. Halmos, and A.~L. Shields.
\newblock Ces\`aro operators.
\newblock {\em Acta Sci. Math. (Szeged)}, 26:125--137, 1965.

\bibitem{BrunoHaslam}
O.~Bruno and M.~Haslam.
\newblock Regularity theory and superalgebraic solvers for wire antenna
  problems.
\newblock {\em SIAM J. Sci. Comput.}, 29(4):1375--1402, 2007.

\bibitem{BrunoKunyansky}
O.~Bruno and L.~Kunyansky.
\newblock A fast, high-order algorithm for the solution of surface scattering
  problems: basic implementation, tests, and applications.
\newblock {\em J. Comput. Phys.}, 169(1):80--110, 2001.

\bibitem{BrunoLintner3}
O.~Bruno and S.~Lintner.
\newblock A high-order integral solver for scalar problems of diffraction by
  screens and apertures in three dimensional space.
\newblock {\em In preparation}.

\bibitem{ChristiansenNedelec}
S.~H. Christiansen and J.-C. N{\'e}d{\'e}lec.
\newblock Preconditioners for the boundary element method in acoustics.
\newblock In {\em Mathematical and numerical aspects of wave propagation
  ({S}antiago de {C}ompostela, 2000)}, pages 776--781. SIAM, Philadelphia, PA,
  2000.

\bibitem{ColtonKress1}
D.~Colton and R.~Kress.
\newblock {\em Integral Equation Methods in Scattering Theory}.
\newblock John Wiley \& Sons, 1983.

\bibitem{ColtonKress2}
D.~Colton and R.~Kress.
\newblock {\em Inverse Acoustic and Electromagnetic Scattering Theory}.
\newblock Springer, 1997.

\bibitem{Bateman}
A.~Erd{\'e}lyi, W.~Magnus, F.~Oberhettinger, and F.~Tricomi.
\newblock {\em Higher transcendental functions. {V}ol. {II}}.
\newblock Robert E. Krieger Publishing Co. Inc., Melbourne, Fla., 1981.
\newblock Based on notes left by Harry Bateman, Reprint of the 1953 original.

\bibitem{StephanHsiao}
George~C. Hsiao, Ernst~P. Stephan, and Wolfgang~L. Wendland.
\newblock On the {D}irichlet problem in elasticity for a domain exterior to an
  arc.
\newblock {\em J. Comput. Appl. Math.}, 34(1):1--19, 1991.

\bibitem{RokhlinJiang}
S.~Jiang and V.~Rokhlin.
\newblock Second kind integral equations for the classical potential theory on
  open surfaces. {II}.
\newblock {\em J. Comput. Phys.}, 195(1):1--16, 2004.

\bibitem{Kress}
R.~Kress.
\newblock {\em Linear integral equations}, volume~82 of {\em Applied
  Mathematical Sciences}.
\newblock Springer-Verlag, New York, second edition, 1999.

\bibitem{Kussmaul}
R.~Kussmaul.
\newblock Ein numerisches {V}erfahren zur {L}\"osung des {N}eumannschen
  {N}eumannschen {A}ussenraumproblems f\"ur die {H}elmholtzsche
  {S}chwingungsgleichung.
\newblock {\em Computing (Arch. Elektron. Rechnen)}, 4:246--273, 1969.

\bibitem{LintnerBruno}
S.~Lintner and O.~Bruno.
\newblock A generalized {C}alder\'on formula for open-arc diffraction problems:
  theoretical considerations.
\newblock {\em Submitted}.

\bibitem{Costabel}
R.~Duduchava M.~Costabel, M.~Dauge.
\newblock Asymptotics without logarithmic terms for crack problems.
\newblock {\em Communications in PDE}, pages 869--926, 2003.

\bibitem{Martensen}
Erich Martensen.
\newblock \"{U}ber eine {M}ethode zum r\"aumlichen {N}eumannschen {P}roblem mit
  einer {A}nwendung f\"ur torusartige {B}erandungen.
\newblock {\em Acta Math.}, 109:75--135, 1963.

\bibitem{MasonHandscomb}
J.~C. Mason and D.~C. Handscomb.
\newblock {\em Chebyshev polynomials}.
\newblock Chapman \& Hall/CRC, Boca Raton, FL, 2003.

\bibitem{Monch}
L.~M{\"o}nch.
\newblock On the numerical solution of the direct scattering problem for an
  open sound-hard arc.
\newblock {\em J. Comput. Appl. Math.}, 71(2):343--356, 1996.

\bibitem{Nedelec}
JC. N{\'e}d{\'e}lec.
\newblock {\em Acoustic and electromagnetic equations}, volume 144 of {\em
  Applied Mathematical Sciences}.
\newblock Springer-Verlag, New York, 2001.

\bibitem{PovznerSuharesvki}
A.~Ya. Povzner and I.~V. Suharevski{\u\i}.
\newblock Integral equations of the second kind in problems of diffraction by
  an infinitely thin screen.
\newblock {\em Soviet Physics. Dokl.}, 4:798--801, 1960.

\bibitem{NumericalRecipes}
W.~H. Press, Teukolsky S., W.~Vetterling, and B.~Flannery.
\newblock {\em Numerical Recipes in C}.
\newblock Cambridge University Press, 1992.
\newblock Second Edition.

\bibitem{Rokhlin}
V.~Rokhlin.
\newblock Diagonal forms of translation operators for the {H}elmholtz equation
  in three dimensions.
\newblock {\em Appl. Comput. Harmon. Anal.}, 1(1):82--93, 1993.

\bibitem{saad}
Youcef Saad and Martin~H. Schultz.
\newblock G{MRES}: a generalized minimal residual algorithm for solving
  nonsymmetric linear systems.
\newblock {\em SIAM J. Sci. Statist. Comput.}, 7(3):856--869, 1986.

\bibitem{VainikkoSaranen}
J.~Saranen and G.~Vainikko.
\newblock {\em Periodic integral and pseudodifferential equations with
  numerical approximation}.
\newblock Springer Monographs in Mathematics. Springer-Verlag, Berlin, 2002.

\bibitem{Stephan}
E.~Stephan.
\newblock Boundary integral equations for screen problems in {${\bf R}\sp 3$}.
\newblock {\em Integral Equations Operator Theory}, 10(2):236--257, 1987.

\bibitem{StephanTran}
E.~Stephan and T.~Tran.
\newblock Domain decomposition algorithms for indefininte hypersingular
  integral equations: the h and p versions.
\newblock {\em SIAM J. Sci. Comp.}, 19(4):1139--1153, 1998.

\bibitem{StephanWendland}
E.~Stephan and W.~Wendland.
\newblock An augmented {G}alerkin procedure for the boundary integral method
  applied to two-dimensional screen and crack problems.
\newblock {\em Applicable Anal.}, 18(3):183--219, 1984.

\bibitem{StephanWendland2}
W.~L. Wendland and E.~P. Stephan.
\newblock A hypersingular boundary integral method for two-dimensional screen
  and crack problems.
\newblock {\em Arch. Rational Mech. Anal.}, 112(4):363--390, 1990.

\bibitem{YanSloan}
Y.~Yan and I.~H. Sloan.
\newblock On integral equations of the first kind with logarithmic kernels.
\newblock {\em J. Integral Equations Appl.}, 1(4):549--579, 1988.

\end{thebibliography}

\end{document}